\documentclass{article}  \usepackage{amssymb,latexsym}
\usepackage{array}
\usepackage[all,2cell,ps]{xy}
\newcommand{\edge}[1]{\ar@{-}[#1]}

\begin{document}
\title{Equivariant K-theory of compactifications of algebraic groups}
\author{V.Uma} \date{} \maketitle \thispagestyle{empty}


\def\theequation {\arabic{section}.\arabic{equation}}
\newcommand{\codim}{\mbox{{\rm codim}$\,$}}
\newcommand{\stab}{\mbox{{\rm stab}$\,$}}
\newcommand{\lr}{\mbox{$\longrightarrow$}}
\newcommand{\ch}{{\cal H}}
\newcommand{\cf}{{\cal F}}
\newcommand{\cd}{{\cal D}}
\newcommand{\blr}{\longrightarrow}
\newcommand{\da}{\Big\downarrow}
\newcommand{\ra}{ \rightarrow}
\newcommand{\ua}{\Big\uparrow}
\newcommand{\hra}{\hookrightarrow}
\newcommand{\bhra}{\mbox{\LARGE{$\hookrightarrow$}}}
\newcommand{\rt}{\mbox{\Large{$\rightarrowtail$}}}
\newcommand{\dua}{\begin{array}[t]{c}
\Big\uparrow \\ [-4mm]
\scriptscriptstyle \wedge \end{array}}
\newcommand{\be}{\begin{equation}}
\newcommand{\ee}{\end{equation}}

\newtheorem{guess}{Theorem}[section]
\newcommand{\bth}{\begin{guess}$\!\!\!${\bf .}~}
\newcommand{\eeth}{\end{guess}}
\renewcommand{\bar}{\overline}
\newtheorem{propo}[guess]{Proposition}
\newcommand{\bpropo}{\begin{propo}$\!\!\!${\bf .}~}
\newcommand{\epropo}{\end{propo}}

\newtheorem{lema}[guess]{Lemma}
\newcommand{\blem}{\begin{lema}$\!\!\!${\bf .}~}
\newcommand{\elem}{\end{lema}}

\newtheorem{defe}[guess]{Definition}
\newcommand{\bdefe}{\begin{defe}$\!\!\!${\bf .}~}
\newcommand{\edefe}{\end{defe}}

\newtheorem{note}[guess]{\bf Notation}

\newtheorem{coro}[guess]{Corollary}
\newcommand{\bcor}{\begin{coro}$\!\!\!${\bf .}~}
\newcommand{\ecor}{\end{coro}}

\newtheorem{rema}[guess]{Remark}
\newcommand{\brem}{\begin{rema}$\!\!\!${\bf .}~\rm}
\newcommand{\erem}{\end{rema}}

\newtheorem{exam}[guess]{Example}
\newcommand{\beg}{\begin{exam}$\!\!\!${\bf .}~\rm}
\newcommand{\eeg}{\end{exam}}

\newcommand{\ctext}[1]{\makebox(0,0){#1}}
\setlength{\unitlength}{0.1mm}
\newcommand{\cl}{{\cal L}}
\newcommand{\cp}{{\cal P}}
\newcommand{\bz}{\mathbb{Z}}
\newcommand{\cs}{{\cal s}}
\newcommand{\cv}{{\cal V}}
\newcommand{\ce}{{\cal E}}
\newcommand{\ck}{{\cal K}}
\newcommand{\cR}{{\cal R}}
\newcommand{\cz}{{\cal Z}}
\newcommand{\cj}{{\cal J}}
\newcommand{\cg}{{\cal G}}
\newcommand{\ci}{{\cal I}}
\newcommand{\bq}{\mathbb{Q}}
\newcommand{\ba}{\mathbb{A}}
\newcommand{\bt}{\mathbb{T}}
\newcommand{\bg}{\mathbb{G}}
\newcommand{\bh}{\mathbb{H}}
\newcommand{\br}{\mathbb{R}}
\newcommand{\tG}{\widetilde{G}}
\newcommand{\tB}{\widetilde{B}}
\newcommand{\tC}{\widetilde{C}}
\newcommand{\tT}{\widetilde{T}}
\newcommand{\tP}{\widetilde{P}}
\newcommand{\wt}{\widetilde}
\newcommand{\tW}{\widetilde{W}}
\newcommand{\tp}{\widetilde{p}}
\newcommand{\tphi}{\widetilde{\Phi}}
\newcommand{\im}{{\rm Im}\,}
\newcommand{\bc}{\mathbb{C}}
\newcommand{\bp}{\mathbb{P}}
\newcommand{\spin}{{\rm Spin}\,}
\newcommand{\ds}{\displaystyle}
\newcommand{\tor}{{\rm Tor}\,}
\newcommand{\bff}{{\bf F}}
\newcommand{\bs}{\mathbb{S}}
\def\ns{\mathop{\lr}}
\def\nssup{\mathop{\lr\,sup}}
\def\nsinf{\mathop{\lr\,inf}}
\renewcommand{\phi}{\varphi}
\newcommand{\co}{{\cal O}}
\newcommand{\cc}{{\cal C}}
\newcommand{\ly}{{L^{Y}}}
\newcommand{\dsi}{{\Delta\setminus I}}
\newcommand{\dsj}{{\Delta\setminus J}}

\begin{abstract}
In this article we describe the $G\times G$-equivariant
$K$-ring of $X$, where $X$ is a regular compactification of a
connected complex reductive algebraic group $G$. Furthermore, in the
case when $G$ is a semisimple group of adjoint type, and $X$ its
wonderful compactification, we describe its ordinary $K$-ring
$K(X)$. More precisely, we prove that $K(X)$ is a free module over
$K(G/B)$ of rank the cardinality of the Weyl group. We further give an
explicit basis of $K(X)$ over $K(G/B)$, and also determine the
structure constants with respect to this basis.
\end{abstract}

\section*{Introduction}

Let $G$ denote a connected complex reductive algebraic group,
$B\subset G$ a Borel subgroup and $T\subset B$ a maximal torus of
dimension $l$. Let $C$ be the center of $G$ and let $G_{ad}:=G/C$ be
the corresponding semisimple adjoint group. Let $W$ denote the Weyl
group of $(G,T)$.

A normal complete variety $X$ is called an {\it equivariant
compactification} of $G$ if $X$ contains $G$ as an open subvariety and
the action of $G\times G$ on $G$ by left and right multiplication
extends to $X$.  We say that $X$ is a {\it regular compactification}
of $G$ if $X$ is an equivariant compactification of $G$ which is
regular as a $G\times G$-variety ( \cite[Section 2.1]{Br}).  Smooth
complete toric varieties are regular compactifications of the
torus. For the adjoint group $G_{ad}$, the wonderful compactification
$\bar {G_{ad}}$ constructed by De Concini and Procesi in \cite{DP1} is
the unique regular compactification of $G_{ad}$ with a unique closed
$G_{ad}\times G_{ad}$-orbit.

The main aim of this article is to describe the $T\times
T$-equivariant and $G\times G$-equivariant $K$-ring of $X$. For this
purpose, we essentially follow the methods used in the description of
the $T\times T$-equivariant and $G\times G$-equivariant Chow ring of
$X$ by Brion (\cite[Section 3]{Br}).  Indeed, we see that these
methods can be naturally generalised to the setting of $K$-theory, for
the purpose of which we apply as key tools, the localisation theorem
of Vezzosi and Vistoli (\cite[Theorem 2]{VV}) and the results of
Merkurjev (\cite[Theorem 4.2]{Mer}).

We begin with a Preliminary section \S1, where we recall basic notions
on equivariant $K$-theory and prove certain necessary facts which are
later used in proving the main results. We refer to \S1 and \S2 for
the notations used below.

In \S2 (see Theorem \ref{kring}) we describe $K_{T\times T}(X)$ in
terms of closed $G\times G$-orbits and the $T\times T$-invariant
curves described in \cite[Section 3]{Br}. More precisely, using the
localisation theorem we embed $K_{T\times T}(X)$ inside $\prod_{\sigma
\in\cf_{+}(l)} K_{T\times T}(Z_{\sigma})$, where each
$Z_{\sigma}\simeq G/B^{-}\times G/B$ is a closed $G\times G$-orbit in
$X$. The image of $K_{T\times T}(X)$ inside
$\prod_{\sigma\in\cf_{+}(l)} K_{T\times T}(Z_{\sigma})$ is further
described by certain equivalence relations which are completely
determined by the $T\times T$-action on the $T\times T$-invariant
curves joining the $T\times T$-fixed points, which are the base points
of the closed orbits.

Using the above, we further get a description of $K_{G\times G}(X)$ in
Cor. \ref{co1} and Cor. \ref{co2}. In particular, we prove in
Cor. \ref{co2} that $K_{G\times G}(X)\simeq (K_{T\times
T}(X))^{W}\simeq (K_{T}(\bar{T})\otimes R(T))^W$, where $\bar{T}$
denotes the closure of $T$ in $X$. As a consequence, $K_{G\times
G}(X)$ is a module over its subring $R(T)\otimes R(T)^{W}\simeq
R(T)\otimes R(G)$. Here we mention that Cor. \ref{co2} is analogous to
the corresponding result for equivariant cohomology of wonderful
compactifications due to Littelmann and Procesi (\cite{LP}).

In Theorem \ref{kdec} we give an explicit description of the additive
structure of $K_{G\times G}(X)$ as a module over its subring $1\otimes
R(G)$.  More precisely, we give a direct sum decomposition of
$K_{G\times G}(X)$, where each piece of the decomposition is a
$1\otimes R(G)$-submodule of the ring $K_{T}(\bar{T}^{+})\otimes R(T)$
(see \S2 for the definition of the toric variety $\bar{T}^{+}$).
Further, by defining the multiplication of the pieces inside the
subring $K_{T}(\bar{T}^{+})\otimes R(T)$ we describe the ring
structure of $K_{G\times G}(X)$, and obtain the explicit
multiplication rule (see Cor. \ref{norm}).  Moreover, from the direct
sum decomposition we also get a natural multifiltration (see
Cor. \ref{kfilt}) where the filtered pieces are $R(T)\otimes
R(G)$-submodules.

The rational equivariant cohomology of regular embeddings of symmetric
spaces have been described by Bifet, De Concini and Procesi
(\cite{BDP}) in terms of Stanley-Reisner systems. Our approach via the
localisation theorem yields another proof of their results for group
embeddings, and also an integral version via $K$-theory.

In \S3, we take $G$ to be the simply connected cover of the semisimple
adjoint group $G_{ad}$, and $T$ a maximal torus of $G$. Then {\it for
the wonderful compactification $X$} of $G_{ad}$, we give a direct sum
decomposition of $K_{G\times G}(X)$ as a free module of rank $|W|$
over the subring $R(T)\otimes R(G)$ (see Theorem \ref{kwond}).
Moreover, each piece of the direct sum is canonically isomorphic to
submodules of $R(T)\otimes R(T)$. This enables us to describe the
multiplication of the direct sum pieces inside the subring
$R(T)\otimes R(T)$.  We also give an explicit description of the
multiplicative structure and the multiplication rule of the basis
elements (see Theorem \ref{kdec1}).

Finally, by further application of the result of Merkurjev
(\cite[Theorem 4.2]{Mer}), we describe the ordinary $K$-ring of
$X$. More precisely, we prove that the subring generated by $Pic(X)$
in $K(X)$ is canonically isomorphic to $K(G/B)$, and $K(X)$ is a free
module of rank $|W|$ over this subring. Furthermore, we also give a
precise description of the multiplication of the basis elements over
$K(G/B)$ in ordinary $K$-ring by pushing down the multiplicative
structure in the equivariant $K$-ring. More precisely, in Theorem
\ref{main} we construct an explicit basis of $K(X)$ over $K(G/B)$ and
determine the structure constants with respect to this basis.

\noindent
{\bf Acknowledgements:} This work was carried out during my one year
stay at Institut Fourier, Grenoble under the French Government
post-doctoral fellowship.  I am grateful to Michel Brion for his
insightful guidance and for the invaluable discussions we have had
throughout the period of this work at Institut Fourier. I also thank
him for a careful reading of earlier versions of this manuscript and
giving several suggestions for improving its presentation.  I thank
V. Balaji for valuable discussions and for constant encouragement
during this work. I thank both the referees for their careful reading
of the manuscript and for their comments and valuable suggestions
which has led to improving the exposition.

\section{Preliminaries}

\subsection{Regular group compactifications}
Let $W$ denote the Weyl group and $\Phi$ denote the root system of
$(G,T)$. We have the subset $\Phi^{+}$ of positive roots and its
subset $\Delta=\{\alpha_1,\ldots, \alpha_r\}$ of simple roots where $r$
is the semisimple rank of $G$. For $\alpha\in\Delta$ we denote by
$s_{\alpha}$ the corresponding simple reflection. For any subset
$I\subset \Delta$, let $W_{I}$ denote the subgroup of $W$ generated by
all $s_{\alpha}$ for $\alpha\in I$. At the extremes we have
$W_{\emptyset}=\{1\}$ and $W_{\Delta}=W$.

A $G$-{\it variety} is a complex algebraic variety with an algebraic
action of $G$.

We now recall the definition of a regular $G$-variety due to Bifet, De
Concini and Procesi. (see \S3 of \cite{BDP} and \S1.4 of \cite{Br}).

\bdefe A $G$-variety $X$ is said to be {\it regular} if it satisfies
the following conditions:

(i) $X$ is smooth and contains a dense $G$-orbit $X^0_G$ whose
complement is a union of irreducible smooth divisors with normal
crossings (the boundary divisors).

(ii) Any $G$-orbit closure in $X$ is the transversal intersection of
the boundary divisors which contain it.

(iii) For any $x\in X$, the normal space $T_xX/T_x(Gx)$ contains a
dense orbit of the isotropy group $G_x$.
\edefe

Consider the connected reductive group $G$ as a homogeneous space
under $G\times G$ for the action given by left and right
multiplication: $(g_1,g_2)\gamma = g_1\gamma g_2^{-1}$. Then the
isotropy group of the identity is the diagonal $diag(G)$.

A normal complete variety $X$ is called an equivariant {\it
compactification} of $G$, if $X$ contains $G$ as an open subvariety and
the action of $G\times G$ on $G$ by left and right multiplication
extends to $X$.

We now recall the definition of a regular compactification of $G$ (see
\S2.1 of \cite{Br}).

\bdefe We say that $X$ is a regular compactification of $G$ if $X$ is
a $G\times G$-equivariant compactification of $G$ which is regular as
a $G\times G$-variety.  \edefe

\noindent
{\bf Examples:}
\begin{enumerate}
\item Smooth complete toric varieties are regular compactifications of
the torus.

\item For the adjoint group $G_{ad}$, the wonderful compactification
$\bar G_{ad}$ constructed by De Concini and Procesi in \cite{DP1} is
the unique regular compactification of $G_{ad}$ with a unique closed
$G_{ad}\times G_{ad}$-orbit.

\end{enumerate}

\subsection{Preliminaries on K-theory}
Let $X$ be a smooth projective complex $G$-variety.  Let $K_{G}(X)$
and $K_{T}(X)$ denote the Grothendieck groups of $G$ and
$T$-equivariant coherent sheaves on $X$ respectively. Recall that
$K_{T}(pt)=R(T)$ and $K_{G}(pt)=R(G)$ where $R(T)$ and $R(G)$ denote
respectively the Grothendieck group of complex representations of $T$
and $G$. The Grothendieck group of equivariant coherent sheaves can be
identified with the Grothendieck ring of equivariant vector bundles on
$X$. Further, the structure morphism $X\ra {Spec}~{\bc}$ induces
canonical $R(G)$ and $R(T)$-module structures on $K_{G}(X)$ and
$K_{T}(X)$ respectively (see Prop. 5.1.28 of \cite{CG}) and Example
2.1 of \cite{Mer}).

Let $\Lambda:=X^*(T)$. Then $R(T)$ (the representation ring of the
torus $T$) is isomorphic to the group algebra $\bz[{\Lambda}]$. Let
$e^{\lambda}$ denote the element of $\bz[{\Lambda}]=R(T)$
corresponding to a weight $\lambda\in {\Lambda}$.  Then
$(e^{\lambda})_{\lambda\in {\Lambda}}$ is a basis of the $\bz$ module
$\bz[{\Lambda}]$. Further, since $W$ acts on $X^*(T)$, on
$\bz[{\Lambda}]$ we have the following natural action of $W$ given by
: $w(e^{\lambda})=e^{w(\lambda)}$ for each $w\in W$ and $\lambda\in
{\Lambda}$. Recall that we can identify $R(G)$ with $R(T)^{W}$ via
restriction to $T$, where $R(T)^W$ denotes the subring of $R(T)$
invariant under the action of $W$.

The following is a theorem analogous to Theorem 3.4 of \cite{Br2},
which we shall use to prove the main result.

\bth{\label{im}} Let $X$ be a nonsingular projective variety on which
$T$ acts with finitely many fixed points $x_1,\ldots,x_m$ and finitely
many invariant curves. Then the image of

$$\iota^{*}: K_{T}(X) \ra K_{T}(X^{T})$$

\noindent
is the set of all $(f_1,\ldots,f_m)\in R(T)^{m}$ such that $f_i\equiv f_j
\pmod {(1-e^{-\chi})}$ whenever $x_i$ and $x_j$ lie in an invariant
irreducible curve $C$ and  $T$-acts on $C$ through the character $\chi$.
\eeth

{\it Proof:} By Theorem 2 of \cite{VV} it follows that the above
restriction homomorphism $\iota^{*}$ is injective and its image is equal
to the intersection of all the images of the restriction homomorphisms
$K_{T}(X^{T^{\prime}})\ra K_{T}(X^{T})$ for all subtori
$T^{\prime}\subseteq T$ of codimension $1$.

Since $X$ contains finitely many invariant curves, $X^{T^\prime}$ is
at most one dimensional for every codimension $1$ subtorus
$T^{\prime}\subset T$.

Let $X_{n-1}:=\bigcup X^{T^{\prime}}$, where the union runs over all
subtori $T^{\prime}$ of codimension one in $T$. Since $X_{n-1}$ is
one-dimensional it consists of disjoint union of points and
nonsingular irreducible curves; let $C$ be such a curve.

If $C$ contains a unique fixed point $x$, then $\iota_x^{*}:K_{T}(C)\ra
K_{T}(x)=R(T)$ is an isomorphism. Otherwise, $C$ is isomorphic to
$\bp^1$. It follows that $C$ contains two distinct fixed points $x$
and $y$. Moreover, the image of

$$\iota^{*}_C : K_{{T}}(C) \ra K_{T}(x)\times K_{T}(y)
=R(T)\times R(T)$$

\noindent
consists of pairs of elements $(f,g)\in R(T)\times R(T)$ such that
$f\equiv g \pmod {(1-e^{-\chi})}$ where $T$ acts on $C$ through the weight
$\chi$. This can be seen as follows:

Let us choose as a basis of $K_{T}(C)$ over $R(T)$ the class of
the trivial line bundle ${{\co}_{C}}$, and the class of the Hopf
bundle $H$ which is the dual of the tautological bundle. Then under
$\iota^{*}_C$, the image of $\co_{C}$ is $(1,1)$ and that of $H$ is
$(e^{\chi},1)$. Since any element in $K_{T}(C)$ is a linear
combination of $\co_{C}$ and $H$, the difference of the coordinates of
the image in $K_{T}(C^{T})$ is always divisible by
$1-e^{-\chi}$. Conversely, let $(f,g)\in R(T)\times R(T)$ be such
that $(1-e^{-\chi})\cdot h= f-g$ for some $h\in R(T)$. Then we see
that the element $g[{\co}_{C}]+ e^{-\chi}\cdot h ([H]-[\co_{C}])$ in
$K_{T}(C)$ maps to $(f,g)$ under $\iota^{*}_{C}$.

The theorem now follows by applying Theorem 2 of \cite{VV}. \hfill $\Box$

Recall from Cor. 3.7 of \cite{Iv} that there exists an exact sequence:
$$1\ra \cz\ra \tG:=\tC \times G^{ss}{\stackrel{\pi}{\lr}}G \ra 1
\eqno(1.1)$$ where $\cz$ is a finite central subgroup, $\tC$ is a
torus and $G^{ss}$ is semisimple and simply-connected.  The
condition that $G^{ss}$ is simply connected implies that $\tG$ is
{\it factorial} (see \cite{Mer}).

Then $\tB:=\pi^{-1}(B)$ and $\tT:=\pi^{-1}(T)$ are respectively a
Borel subgroup and a maximal torus of $\tG$. Further, by restricting
the map $\pi$ to $\tT$ we get the following exact sequence:
$$ 1\ra \cz\ra \tT\ra T \ra 1. \eqno(1.2)$$

Let $\tW$ and $\tphi$ denote respectively the Weyl group and the root
system of $(\tG,\tT)$.  Then by the exact sequence $(1.1)$ it follows
in particular that $\tW=W$ and $\tphi=\Phi$.

Further we have $R(\tG)= R(\tC)\otimes R(G^{ss})$ and $R(\tT)\simeq
R(\tC)\otimes R(T^{ss})$ where $T^{ss}$ is the (unique) maximal torus
$\tT\cap G^{ss}$.

Recall we can identify $R(\tG)$ with $R(\tT)^{W}$ via restriction to
$\tT$, and further $R(\tT)$ is a free $R(\tG)$ module of rank $|W|$
(see Theorem 6.41, pp.164 of \cite{Ad} and Theorem 1, pp. 199 of
\cite{Pit}). Moreover, since $G^{ss}$ is semi-simple and simply
connected, $R(G^{ss})\simeq \bz[x_1,\ldots, x_r]$ is a polynomial
ring on the ¨fundamental representations¨. Hence $R(\tG)=R(\tC)\otimes
R(G^{ss})$ is the tensor product of a polynomial ring and a Laurent
polynomial ring, and hence a regular ring of dimension
$r+{dim}(\tC)={rank}(G)$ where $r$ is the rank of $G^{ss}$.

We shall consider the $\tT$ and $\tG$-equivariant $K$-theory of $X$
where we take the natural actions of $\tT$ and $\tG$ on $X$ through
the canonical surjections to $T$ and $G$ respectively.

We consider $\bz$ as an $R(\tG)$-module by the augmentation map
$\epsilon:R(\tG)\ra \bz$ which maps any $\tG$-representation $V$ to
${dim}(V)$. Moreover, we have the natural restriction
homomorphisms $K_{\tG}(X)\ra$ $K_{\tT}(X)$ and $K_{\tG}(X)\ra K(X)$
where $K(X)$ denotes the ordinary Grothendieck ring of algebraic
vector bundles on $X$.  We then have the following isomorphisms (see
Prop. 4.1 and Theorem 4.2 of \cite{Mer}) (also see Theorem 6.1.22
pp.310 of \cite{CG}):

\begin{enumerate}
\item[$(a)$] $R(\tT)\otimes_{R(\tG)}K_{\tG}(X)\simeq K_{\tT}(X).$

\item[$(b)$] $K_{\tG}(X)\simeq K_{\tT}(X)^{W}$.

\item[$(c)$] $\bz\otimes_{R(\tG)} K_{\tG}(X)\simeq K(X).$
\end{enumerate}

\brem{\label{hk}}
In fact the above isomorphisms $(a)$ and $(b)$ hold in higher
equivariant K-theory and the isomorphism $(c)$ corresponds to the
degeneration of the Merkurjev spectral sequence
$E^{2}_{p,q}=Tor_p^{R(\tG)}(\bz, K^{q}_{\tG}(X))\Rightarrow
K_{p+q}(X)$ (see pp. 2-3 of \cite{Mer}).
\erem

\brem We will prove in Theorem \ref{inv} that the isomorphism $(b)$
also holds when $\tG$ and $\tT$ are replaced with $G$ and $T$
respectively i.e, $K_{G}(X)\simeq K_{T}(X)^{W}$.  \erem

\blem\label{fm} Let $X$ be a smooth projective $G$-variety containing
only finitely many $T$-fixed points. Then $K_{T}(X)$ is a free module
over $R(T)$ of rank $|X^{T}|$. Furthermore, $K_{\tT}(X)$
(resp. $K_{\tG}(X)$) is also free over $R(\tT)$ (resp. $R(\tG)$) of
rank $|X^{T}|$.  \elem

\noindent {\it Proof:} Since $X$ is a smooth projective variety with
$T$-action such that $X^{T}$ is finite, it admits a Bialynicki-Birula
cellular decomposition with $m=|X^{T}|$ $T$-stable affine cells.  Let
$X_1, X_2,\ldots, X_m=pt$ be an ordering of the cells with
${dim}(X_1)\geq {dim}(X_2)\geq \ldots$. Set $X^{j}=\sqcup_{i\geq j}
X_i$. Then $X=X^0\supset X^1\supset\cdots\supset X^m=pt$ is a
decreasing filtration on $X$ by closed $T$-stable subvarieties. Thus
$X\ra pt$ is a $T$-equivariant cellular fibration over a
point. Therefore by the Cellular Fibration Lemma (see pp. 270
\cite{CG}) it follows that, $K_{T}(X)$ is free module over $K_{T}(pt)=
R(T)$ of rank $m$.

Since $\tT$ acts on $X$ via the canonical surjection to $T$ it
similarly follows that $K_{\tT}(X)$ is free module over $K_{\tT}(pt)=
R(\tT)$ of rank $|X^{\tT}|=|X^{T}|=m$.

Now, since $K_{\tT}(X)$ is a free module over $R(\tT)$, and $R(\tT)$
is free over $R(\tG)$, it follows that $K_{\tT}(X)$ is a free module
over $R(\tG)$.  Further, since $R(\tG)$ is a direct summand of
$R(\tT)$, the isomorphism $(a)$ above implies that $K_{\tG}(X)$ is a
direct summand of $K_{\tT}(X)$ as an $R(\tG)$-module. Thus
$K_{\tG}(X)$ is a projective module over $R(\tG)$. Moreover, since
$R(\tG)$ is a tensor product of a polynomial ring and a Laurent
polynomial ring $K_{\tG}(X)$ is in fact free over $R(\tG)$ (see
Theorem 1.1 of \cite{Gu}).

The isomorphism $(a)$ above further implies that the rank of
$K_{\tG}(X)$ over $R(\tG)$ is same as the rank of $K_{\tT}(X)$ over
$R(\tT)$ which is $m$.  \hfill $\Box$

\blem{\label{ext}} Let $1\ra \cz\ra\tT\ra T\ra 1$ be an exact sequence
of algebraic groups where $\tT$ and $T$ are complex tori and $\cz$ is
a finite abelian group. Let $1_{\cz}={\chi}_1,\ldots, {\chi}_m$ denote
the characters of $\cz$. Further, let $1_{\tT}=\wt{{\chi}_1},\ldots,
\wt{{\chi}_m}$ be arbitrary lifts of ${\chi}_1,\ldots, {\chi}_m$ to
characters of $\tT$. Let $Y$ be an irreducible $T$-variety. We then
have the following isomorphisms:

\begin{enumerate}

\item[(i)] $R(\tT)\simeq \bigoplus_{i=1}^{m}
R(T)e^{\widetilde{\chi_i}}$

\item[(ii)] $K_{\tT}(Y)\simeq K_T(Y)\otimes_{R(T)} R(\tT)$

\end{enumerate}
\elem

\noindent
{\it Proof:} Let $V$ be any $\tT$-representation and
$V=\bigoplus_{\chi\in X^{*}(\tT)}V_\chi$ be the direct sum
decomposition of $V$ as $\tT$-weight spaces. Then
$$V_{i}:=\bigoplus_{\chi\mid_{\cz}= \chi_{i}} V_{\chi}.$$ are the
isotypical components with respect to the characters $\chi_1,\ldots,
\chi_m$ of $\cz$.  Thus as $\tT$-modules we have an isomorphism
$V_{i}\simeq e^{\wt{\chi_i}}\otimes V^{i}$, where the $\tT$-module
$V^i$ is in fact a $T$-module since the $\cz$-action on it is trivial.
Since $V=\bigoplus_{i=1}^m V_i$, it follows that as $\tT$-modules we
have $V\simeq \bigoplus_{i=1}^m e^{\wt{\chi_i}}\otimes V^{i}$. This
proves $(i)$.

We have a canonical homomorphism of rings $R(\tT)\otimes_{R(T)}
K_{T}(Y)\ra K_{\tT}(Y)$ where, $[V]\in R(\tT)$ maps to the trivial
bundle $Y\times V$ and the map from $K_{T}(Y)$ to $K_{\tT}(Y)$ is
induced by the surjection $\tT\ra T$.  To define the inverse of the
above homomorphism:

Let $E$ be a $\tT$-equivariant vector bundle on $Y$. Since $Y$ is a
$T$-variety, the $\cz$-action on $Y$ is trivial. Thus on every fibre of
$E$ we get a canonical linear $\cz$-action, which gives a weight space
decomposition on each fiber. Note that since $\cz$ is finite, the weights
of $\cz$ form a finite set. Moreover, since $E$ is locally trivial the
$\cz$-representation is locally constant and hence globally constant
over the irreducible base $Y$.  Thus we get the following vector
bundle direct sum decomposition
$$E=\bigoplus_{i=1}^m E_i$$ where $E_{i}$ denotes the subbundle whose
fibre is the eigenspace corresponding to the character $\chi_{i}$ of
$\cz$. Thus as $\tT$-equivariant bundles we have an isomorphism
$E_i\simeq e^{\wt{\chi_i}}\otimes E^{i}$, where the $\tT$-equivariant
bundle $E^{i}$ is in fact a $T$-equivariant bundle since the
$\cz$-action on it is trivial.

Therefore the inverse map is defined by sending $E$ to the element
$\bigoplus_{i=1}^m e^{\wt{\chi_i}}\otimes E^{i}$ of
$R(\tT)\otimes_{R(T)} K_{T}(Y)$.  This proves $(ii)$.

\bth{\label{inv}} The restriction homomorphism $K_{G}(X)\ra K_T(X)$
induces an isomorphism $K_{G}(X)\simeq K_{T}(X)^{W}$ where
$K_{T}(X)^{W}$ denotes the subring of $W$-invariants of $K_{T}(X)$.
(For the corresponding result in topological $K$-theory see
\cite[Theorem 4.4]{McL}).

\eeth

{\it Proof:} Recall (see Prop. 2.10 of \cite{Mer}) that we have the
following isomorphism: $$K_{T}(X)\simeq K_{G}(X\times G/B).$$

Therefore the projection $p:X\times G/B \ra X$ induces the pull-back
map $p^*:K_{G}(X)\ra K_{T}(X)$. Note that $p$ is a proper map since
its fibre $G/B$ is complete. Thus we further have the push-forward map
$p_*:K_{T}(X)\ra K_{G}(X)$. Now by the projection formula we get
$p_{*}\circ p^{*}=id$ (see 5.2.13 and 5.3.12 of \cite{CG}). In
particular, it follows that $p^{*}$ is injective.

Similarly, the projection $\tp:X\times \tG/\tB \ra X$ induces the
pull-back and push-forward maps: $\tp^*:K_{\tG}(X)\ra K_{\tT}(X)$ and
$\tp_{*}:K_{\tT}(X)\ra K_{\tG}(X)$ respectively. Further, by the
projection formula we get $\tp_{*}\circ \tp^{*}=id$, and hence
${\tp}^{*}$ is injective. Furthermore, by the isomorphism $(b)$ above
we know that the image of $K_{\tG}(X)$ under $\tp^{*}$ is
$K_{\tT}(X)^{W}$.

Let $u:=\pi:\tG \ra G$ and $v:=\pi\mid_{\tT}:\tT \ra T$. Then $u$ and
$v$ induce the ring homomorphisms: $u^{*}:K_{G}(X)\ra K_{\tG}(X)$ and
$v^{*}:K_{T}(X)\ra K_{\tT}(X)$ respectively. Further, since the
isomorphism $K_{T}(X)\simeq K_{G}(X\times G/B)$ is canonical, $\pi^{*}:
K_{G}(X \times G/B)\ra K_{\tG}(X \times \tG/\tB)$ can be identified
with $v^*$.

Now, by (ii) of Lemma \ref{ext} above, the map $v^*: K_{T}(X)\hra
K_{\tT}(X)$ induced by the surjection $\tT\ra T$ is injective. We now
claim that $v^{*}p^{*}=\tp^{*}u^{*}$, so that $u^*$ is also
injective. This can be seen as follows:

For any $G$-vector bundle $V$ on $X$, let $\wt{V}$ denote $V$ thought
of as a $\tG$-vector bundle via the surjection $\pi:\tG\ra G$.  Then
we see that $p^{*}([V])=[V]\boxtimes [\co_{G/B}]\in K_{G}(X\times
G/B)$.  Further, $u^{*}p^{*}([V])=[\wt{V}]\boxtimes [\co_{\tG/\tB}]\in
K_{\tG}(X\times \tG/\tB)$. Moreover, since $u^{*}([V])=[\wt{V}]$, we
see that $\tp^{*}u^{*}([V])=[\wt{V}]\boxtimes [\co_{\tG/\tB}]\in
K_{\tG}(X\times \tG/\tB)$. Hence the claim.

Using the isomorphism $G/B\simeq \tG/\tB$, and the fact that the push
forward map is functorial it follows that: $\tp_{*}:K_{\tG}(X\times
\tG/\tB)\ra K_{\tG}(X)$ restricts to $p_{*}:K_{G}(X\times G/B)\ra
K_{G}(X)$. That is, $\tp_{*}v^{*}=u^{*}p_{*}$.

Thus we get the following commuting diagram:

\[
\begin{array}{cccccccccc}
K_{G}(X) &{\stackrel {u^{*}}\bhra} &K_{\tG}(X)\\ p_{*}\ua\da p^{*} &
&{\tp_{*}}\ua \da {\tp^{*}}\\ K_{T}(X) &{\stackrel{v^{*}} \bhra}
&K_{\tT}(X)\\
\end{array}
\]

In particular, it follows from the above diagram that
$p^{*}(K_{G}(X))\subseteq K_{T}(X)^{W}$. Hence it remains to show that
$K_{T}(X)^{W}\subseteq p^{*}(K_{G}(X))$. This can be seen as follows:

Let $\alpha\in K_{T}(X)^{W}$. Then $v^{*}\alpha\in
K_{\tT}(X)^{W}$. Further, let $v^{*}\alpha=\tp^{*}\beta$ for $\beta\in
K_{\tG}(X)$. Thus
$u^*p_{*}\alpha=\tp_{*}v^{*}\alpha=\tp_{*}\tp^{*}\beta=\beta$.

Now to show that $\alpha=p^{*}(\gamma)$ for $\gamma\in K_{G}(X)$. If
this is true, this further implies that
$\gamma=p_{*}(\alpha)$. Therefore it is enough to show that
$\alpha=p^{*}p_{*}\alpha$. Since $v^{*}$ is injective this is further
equivalent to showing that $v^*\alpha=v^{*}p^{*}p_{*}\alpha$. But this
follows from the above arguments, since $v^*\alpha=\tp^{*}\beta=
\tp^{*}u^{*}p_{*}\alpha=v^{*}p^{*}p_{*}\alpha$. Hence the
theorem.\hfill $\Box$

Let $R(\tT)^{W_I}$ denote the invariant subring of the ring $R(\tT)$ under
the action of the subgroup $W_I$ of $W$ for every $I\subset
\Delta$. Thus in particular we have, $R(\tT)^{W}=R(\tG)$ and
$R(\tT)^{\{1\}}=R(\tT)$. Further, for every $I\subset \Delta$,
$R(\tT)^{W_I}$ is a free module over $R(\tG)=R(\tT)^{W}$ of rank
$|W/W_{I}|$ (see Theorem 2.2 of \cite {St}). Indeed, Theorem 2.2 of
\cite{St} which we apply here holds for $R(T^{ss})$.  However, since
$W$ acts trivially on the central torus $\tC$ and hence trivially on
$R(\tC)$ we have $R(\tT)^{W_I}=R(\tC)\otimes R(T^{ss})^{W_{I}}$ for
every $I\subseteq \Delta$, and hence we obtain the analogous statement
for $R(\tT)$.

Let $W^{I}$ denote the set of minimal length coset representatives of
the parabolic subgroup $W_I$ for every $I\subset \Delta$. Then
$$W^{I}:=\{w\in W \mid l(wv)=l(w)+l(v) ~\forall~ v\in W_I\}=\{w\in
W\mid w(\Phi_{I}^{+})\subset \Phi^{+}\}$$ where $\Phi_I$ is the root
system associated to $W_I$, with $I$ as the set of simple roots.
Recall (see pp. 19 of \cite{hum}) that we also have:
$$W^{I}=\{w\in W\mid l(ws)>l(w)~for~all~s\in I\}.$$

Note that $J\subseteq I$ implies that $W^{\Delta\setminus
J}\subseteq W^{\Delta\setminus I}$.  Let
$$C^{I}:=W^{\Delta\setminus I}\setminus (\bigcup_{J\subsetneq
I}W^{\Delta \setminus J}).\eqno(1.3)$$

Let $\alpha_1,\ldots,\alpha_r$ be an ordering of the set $\Delta$ of
simple roots and $\omega_1,\ldots,\omega_r$ denote the corresponding
fundamental weights for the root system of $(G^{ss},T^{ss})$.
Since $G^{ss}$ is simply connected, the fundamental weights form a
basis for $X^*(T^{ss})$ and hence for every $\lambda\in
X^*(T^{ss})$, $e^{\lambda}\in R(T^{ss})$ is a monomial in the
elements $e^{\omega_i}:1\leq i\leq r$.

In Theorem 2.2 of \cite{St} Steinberg has defined a basis
$\{f_{v}^{^I}: v\in W^{I}\}$ of $R(T^{ss})^{W_{I}}$ as an
$R(T^{ss})^W$-module.  We recall here this definition: For $v\in
W^{I}$ let
$$p_{v}:=\prod_{v^{-1}\alpha_{i}<0} e^{\omega_{i}}\in
R(\tT).\eqno(1.4)$$ Then
$$f_v^{^{I}}:=\sum_{x\in W_{I}(v)\big{\backslash} W_{I}}
x^{-1}v^{-1}p_{v}\eqno(1.5)$$ where $W_{I}(v)$ denotes the stabilizer
of $v^{-1}p_v$ in $W_{I}$.

We shall also denote by $\{f_v^{^I}:v \in W^{I}\}$ the corresponding
basis of $R(\tT)^{W_{I}}$ as an $R(\tT)^{W}$-module where it is
understood that $f_v^{^I}:=1\otimes f_v^{^I} \in R(\tC)\otimes
R(T^{ss})^{W_{I}}$.

We now fix the following notations before we state the next proposition.

\begin{enumerate}
\item[(a)] For $v\in W^{I}$, we shall denote by $W_{I}^{^\ell}(v)$ the
{\it minimal length representatives} of the cosets in
$W_{I}(v)\big{\backslash} W_{I}$. Note that each $w\in W_I$ can be
uniquely expressed as $w=ux$ where $x\in W_{I}^{^\ell}(v)$ ($x$
is the unique element of smallest length in the coset $W_{I}(v)w$)
and $u\in W_{I}(v)$, such that $l(w)=l(u)+l(x)$.

\item[(b)] Let $I\supseteq J$. For each $x\in
W_{\dsj}^{^\ell}(v)\subseteq W_{\dsj}$ we can now consider the {\it
minimal length representative} of the coset $xW_{\dsi}\in
W_{\dsj}\big{/}W_{\dsi}$ which we shall denote by $x^{\prime}$.
Let
$$\big{[}W_{\dsj}^{^\ell}(v)\big{]}^{^\dsi}:=\{x^{\prime}~:~x\in
W_{\dsj}^{^\ell}(v)\}.$$

\end{enumerate}


\bpropo{\label{modstein}} With the above notations we have the
following :

\begin{enumerate}

\item For $v\in W^{\Delta\setminus I}$ we have:
$$f_{v}^{^{\Delta\setminus I}}=\sum_{x\in W_{\Delta\setminus
I}^{^\ell}(v)}x^{-1}v^{-1}p_{v}=\sum_{x\in W_{\Delta\setminus
I}^{^\ell}(v)}f_{vx}^{^{\emptyset}}$$ where $f_{vx}^{^{\emptyset}}$ is
well defined since $vx\in W^{\emptyset}=W$.

\item For $v\in W^{\Delta\setminus J}$ and for $I\supseteq J$ we have:
$$f_{v}^{^{\Delta\setminus J}}=\sum_{x^{\prime}\in
\big{[}W_{\dsj}^{^\ell}(v)\big{]}^{^\dsi}}
f_{vx^{\prime}}^{^{\Delta\setminus I}}$$ where
$f_{vx^{\prime}}^{^{\Delta\setminus I}}$ is well defined since
$vx^{\prime}\in W^{\Delta\setminus I}$.

\item For $v\in W^{\Delta\setminus I}$, $f_{v}^{^{\Delta\setminus I}}$
is in the $R(\tT)^{W}$-span of $\{f_{v^{\prime}}^{^{\Delta\setminus
J}}: v^{\prime}\in C^{J}\}_{J\subseteq I}$.
\end{enumerate}

\epropo

{\it Proof:} This proposition may be well known to experts but since
we could not find a proof in the literature we give it below.

Let $v\in W^{\dsi}$ and $x\in W^{^{\ell}}_{\Delta\setminus
I}(v)$. Then we claim that :
$$x^{-1}v^{-1}p_v=x^{-1}v^{-1}p_{vx}.\eqno(1.6)$$ We shall prove (1.6)
by induction on the length of $x$. Let $l(x)=1$. Then $x=s_{\beta}$
for $\beta\in \dsi$ such that $s_{\beta}\notin W_{\dsi}(v)$. Thus we
require to show that
$$s_{\beta}v^{-1}p_v=(vs_{\beta})^{-1}p_{vs_{\beta}}.\eqno(1.7)$$ This
is equivalent to showing that $p_{v}=p_{vs_{\beta}}$ which can be seen
as follows:

For $w\in W$ let $R(w):=\{\alpha\in \Phi^{+}: w(\alpha)<0\}$. Then by
$(*)$ in pp. 407 of \cite{jos} it follows that:
$$R(s_{\beta}\cdot v^{-1})=R(v^{-1})\bigsqcup
vR(s_{\beta}).\eqno(1.8)$$ Note that $R(s_{\beta})=\beta$. Moreover,
observe that $v(\beta)$ is not a simple root (for this see below
\footnote{\sf If $v(\beta)$ were a simple root then by (1.8) we have
$v(\beta)\notin R(v^{-1})$. Thus by (1.4) we have
$s_{v(\beta)}p_{v}=p_v$ and hence $vs_{\beta}v^{-1}p_{v}=p_{v}$. This
implies that $s_{\beta}v^{-1}p_{v}=v^{-1}p_v$ which is a contradiction
to our assumption that $s_{\beta}\notin W_{\dsi}(v)$.}).  Hence by
(1.8) it follows that the simple roots in $R(s_{\beta}v^{-1})$ = the
simple roots in $R(v^{-1})$. This by (1.4) further implies that
$p_{v}=p_{vs_{\beta}}$ which proves (1.7).


Now we assume by induction that (1.6) holds for all $y\in
W^{^\ell}_{\dsi}(v)$ with $l(y)<l(x)$.

Let $x=ys_{\beta}$ be a reduced expression for $x$ where {\em $\beta$
is a simple root}, $y\in W_{\dsi}$ and $l(y)=l(x)-1$.  Further, since
$x\in W^{^\ell}_{\dsi}(v)$ we must have $y\notin W_{\dsi}(v)$. Indeed
it can be seen that $y\in W^{^\ell}_{\dsi}(v)$ (for this see below
\footnote{\sf Suppose $y\notin W^{^\ell}_{\dsi}(v)$ then we can
uniquely express $y=zy_{o}$ for $y_{o}\in W^{^\ell}_{\dsi}(v)$ and
$z\in W_{\dsi}(v)$ such that $l(y)=l(z)+l(y_{o})$. Let
$x_{o}=y_{o}s_{\beta}$ so that $x=zx_{o}$. Now if $l(y_{o})\lneq
l(y)$ then $l(x_{o})\leq l(y_{o})+1 \lneq l(y)+1=l(x)$. This
contradicts that $x\in W^{^\ell}_{\dsi}(v)$. Thus it follows that
$l(y)=l(y_{o})$ and hence $y=y_{o}$.}).

Hence by induction assumption
$$y^{-1}v^{-1}p_v=y^{-1}v^{-1}p_{vy}.$$ Let $\Delta\setminus
I_{1}:=\{\alpha\in\dsi: l(ys_{\alpha})>l(y)\}$.  Then we {\it claim}
that:
\begin{enumerate}
\item[(i)] $vy\in W^{\Delta\setminus {I_{1}}}$
\item[(ii)] $s_{\beta}\in W^{^\ell}_{\Delta\setminus
{I_1}}(vy)$
\end{enumerate}

Once we prove the above {\it claim} we see that the equality (1.7) for
$v$ replaced by $vy$ and $I$ replaced by $I_1$ will imply the equality
(1.6). Thus it only remains to prove (i) and (ii) above.

{\it Proof of} (i): Since $v\in W^{\dsi}$ we have:
$l(vys_{\alpha})=l(v)+l(ys_{\alpha})>l(v)+l(y)=l(vy)$ for every
$\alpha\in \Delta\setminus I_{1}$. Hence (i) follows.

{\it Proof of} (ii): Since $x=ys_{\beta}$ is a reduced expression,
clearly $s_{\beta}\in W_{\Delta\setminus I_{1}}$. Suppose that
$s_{\beta}\in W_{\Delta\setminus I_{1}}(vy)$. Since by induction we
have $p_{v}=p_{vy}$, it follows that: $ys_{\beta}y^{-1}\in
W_{\Delta\setminus I}(v)$. This further implies that $x=ys_{\beta}=zy$
for an element $z\in W_{\dsi}(v)$. Since $l(y)\lneq l(x)$ this clearly
contradicts that $x\in W^{^\ell}_{\dsi}(v)$. Thus we conclude that
$s_{\beta}\notin W_{\Delta\setminus I_{1}}(vy)$ which implies (ii).

Observe that, on the right hand side of (1.5), without loss of
generality we can assume that $x\in W_{\Delta\setminus I}^{^\ell}(v)$.
Now (1.5) and (1.6) together imply:
$$f_{v}^{^{\Delta\setminus I}}=\sum_{x\in W_{\Delta\setminus
I}^{^\ell}(v)}x^{-1}v^{-1}p_{v}=\sum_{x\in W_{\Delta\setminus
I}^{^\ell}(v)}(vx)^{-1}p_{vx}=\sum_{W_{\Delta\setminus
I}^{^\ell}(v)}f_{vx}^{^{\emptyset}}\eqno(1.9)$$ which proves (1) of
Prop. \ref{modstein}.

Now let $J\subseteq I$. Further, let $v\in W^{\dsj}$ and $x\in
W^{^\ell}_{\dsj}(v)$. Then we can uniquely express $x=x^{\prime}y$
where $y\in W_{\dsi}$ and $x^{\prime}\in
\big{[}W_{\dsj}^{^\ell}(v)\big{]}^{^\dsi}$.

Since there is no reduced expression of $x^{\prime}$ ending in
$s_{\alpha}$ for $\alpha\in \Delta\setminus I$, it follows that
$vx^{\prime}\in W^{\dsi}$ (for this see below \footnote{\sf since
$l(vx^{\prime}s_{\alpha})=l(v)+l(x^{\prime}s_{\alpha})>l(v)+l(x^{\prime})
=l(vx^{\prime})$}).

Suppose $x=x^{\prime}y$ and $x_1=x^{\prime}y_1$ for $x,x_1\in
W^{^\ell}_{\dsj}(v)$ and $y,y_1\in W_{\dsi}$.  Then we see that $x\neq
x_1$ $\Leftrightarrow$ $W_{\dsi}(vx^{\prime})y\neq
W_{\dsi}(vx^{\prime})y_1$ (for this see below \footnote{\sf Since
$l(x)=l(x^{\prime})+l(y)$ and $x\in W^{^\ell}_{\dsj}(v)$, it follows
that $x^{\prime}\in W^{^\ell}_{\dsj}(v)$. Thus by (1.6)
$p_v=p_{vx^{\prime}}$. Hence $x^{-1}v^{-1}p_v=x_1^{-1}v^{-1}p_v$
$\Leftrightarrow$ $y^{-1}(vx^{\prime})^{-1}p_{vx^{\prime}}
=y_1^{-1}(vx^{\prime})^{-1}p_{vx^{\prime}}$.}).

We can now express (1.9) for $J$ as follows:
$$f_{v}^{^\dsj}=\sum_{x^{\prime}\in
\big{[}W_{\dsj}^{^\ell}(v)\big{]}^{^\dsi}}\sum_{y
}y^{-1}(vx^{\prime})^{-1}p_{vx^{\prime}y}=\sum_{x^{\prime}\in
\big{[}W_{\dsj}^{^\ell}(v)\big{]}^{^\dsi}}
f_{vx^{\prime}}^{^\dsi}\eqno(1.10)$$ which proves (2) of
Prop. \ref{modstein}. (In the above summation $y\in
W_{\dsi}(vx^{\prime})\big{\backslash}W_{\dsi}$.)

Let $v\in W^{\dsi}$. Then $v\in C^{J}$ for a unique $J\subseteq I$.
Now, (1.10) can be expressed as:
$$f_{v}^{^\dsj}=f_{v}^{^\dsi}+\sum_{x^{\prime}\neq
1}f_{vx^{\prime}}^{^\dsi}\eqno(1.11)$$ Since $v\in W^{\dsj}$ and
$x^{\prime}\in W_{\dsj}$ we have $l(vx^{\prime})>l(v)$.

We now {\it claim} that if $v_{_1}$ is of maximal length in $W^{\dsi}$
then in fact $v_{_1}\in C^{I}$. This can be seen as follows:

Suppose $v_{1}\notin C^{I}$. Let $v_{_1}\in W^{\dsj}$ for $J\subsetneq
I$ then $l(v_{_1}s_{\alpha})>l(v_{_1})$ for every $s_{\alpha}\in
W_{\dsj}$. In particular, if $\alpha\in (\dsj)\setminus(\dsi)$ then we
note that $v_{_1}s_{\alpha}\in W^{\dsi}$ (for this see below
\footnote{\sf For if $\beta\in\dsi$, then
$l(v_{_1}s_{\alpha}s_{\beta})=l(v_{_1})+l(s_{\alpha}s_{\beta})>l(v_{_1})+
l(s_{\alpha})=l(v_{_1}s_{\alpha})$, since $s_{\beta}\neq
s_{\alpha}$.}).  Since $l(v_{_1}s_{\alpha})>l(v_1)$ this is a
contradiction to the assumption that $v_{_1}$ is of maximal length in
$W^{\dsi}$.

Thus trivially $f^{^\dsi}_{v_{_1}}$ is in the $R(\tT)^{W}$-span of
$\{f^{^\dsj}_{v^{\prime}}:v^{\prime}\in C^{J}\}_{J\subseteq I}$.  Now
by a decreasing induction on $l(v)$ we can therefore assume that if
$l(v_{_1})>l(v)$ then $f^{^\dsi}_{v_{_1}}$ belongs to the
$R(\tT)^{W}$-span of $\{f^{^\dsj}_{v^{\prime}}:v^{\prime}\in
C^{J}\}_{J\subseteq I}$. Thus (1.11) and the induction assumption
together imply (3) of Prop. \ref{modstein}.\hfill $\Box$

\blem{\label{ds}} For $I\subseteq \Delta$, let
$\{f_v^{^{\Delta\setminus I}}:v\in W^{\Delta\setminus I}\}$ denote the
basis defined by Steinberg of $R(\tT)^{W_{\Delta\setminus I}}$ as an
$R(\tT)^{W}$-module. Recall from (1.3) that $W^{\Delta\setminus
I}=\bigsqcup_{J\subseteq I} C^J$. Then $\{f_{v}^{^{\Delta\setminus
J}}: v\in C^{J} \}_{J\subseteq I}$ also form a basis of
$R(\tT)^{W_{\Delta\setminus I}}$ as $R(\tT)^{W}$-module.  Moreover, if
$$R(\tT)_{I}:=\bigoplus_{v\in C^I}R(\tT)^{W}\cdot
f_v^{^{\Delta\setminus I}} \eqno(1.12)$$
then for every $I\subseteq \Delta$ we have the following direct sum
decomposition as $R(\tT)^W$ modules:
$$R(\tT)^{W_{\Delta\setminus I}}=\bigoplus_{J\subseteq I}
R(\tT)_{J}.\eqno(1.13)$$ This further implies that:
$$R(\tT)^{W_{\Delta\setminus I}}=(\sum_{J\subsetneq I}
R(\tT)^{W_{\Delta\setminus J}})\bigoplus R(\tT)_I.\eqno(1.14)$$ \elem

{\it Proof:} By Prop.\ref{modstein} (3) it follows that
$\{f_{v}^{^{\Delta\setminus J}}: v\in C^{J}\}_{J\subseteq I}$ span
$R(\tT)^{W_{\Delta\setminus I}}$ as $R(\tT)^{W}$-module. It is not
hard to see that $\{f_{v}^{^{\Delta\setminus J}}: v\in
C^{J}\}_{J\subseteq I}$ in fact form a basis of
$R(\tT)^{W_{\Delta\setminus I}}$ as $R(\tT)^{W}$-module (for this see
below \footnote{\sf This is because $R(\tG)=R(\tT)^{W}$ is a domain
and $R(\tT)^{W_{\Delta\setminus I}}$ is a free $R(\tT)^{W}$-module of
rank $|W^{\Delta\setminus I}|$, it can be seen that
$\{f_{v}^{^{\Delta\setminus J}}: v\in C^{J} \}_{J\subseteq I}$ are
linearly independent over $R(\tT)^{W}$ and hence form a basis of
$R(\tT)^{W_{\Delta\setminus I}}$ as $R(\tT)^{W}$-module for every
$I\subseteq \Delta$.}).

Since by (1.3) $W^{\Delta\setminus I}=\bigsqcup_{J\subseteq I} C^J$,
we therefore have the following direct sum decomposition:
$$R(\tT)^{W_{\Delta\setminus I}}=\bigoplus_{J\subseteq
I}\bigoplus_{v\in C^J}R(\tT)^{W}\cdot f_v^{^{\Delta\setminus J}} .$$
Hence by (1.12) we further have:
$$R(\tT)^{W_{\Delta\setminus I}}=\bigoplus_{J\subseteq I} R(\tT)_{J}$$
for every $I\subseteq \Delta$. Now it follows by induction that
$$R(\tT)^{W_{\Delta\setminus I}}=(\sum_{J\subsetneq I} R(\tT)^{W_{\Delta
\setminus J}})\bigoplus R(\tT)_{I} .$$ \hfill $\Box$

\brem{\label{compbasis}} In Lemma \ref{ds} we prove that the Steinberg
basis elements for $R(\tT)^{W_{\Delta\setminus J}}$ which correspond
to the indexing set $C^{J}\subseteq W^{\Delta\setminus I}$ for each
$J\subseteq I$, together form another $R(\tT)^{W}$-basis for
$R(\tT)^{W_{\Delta\setminus I}}$. The difference between the Steinberg
basis for $R(\tT)^{W_{\Delta\setminus I}}$ and the new basis is the
following: for all the elements of the Steinberg basis the superscript
$\Delta\setminus I$ remains constant as the index $v$ in the subscript
varies over the elements of $W^{\Delta\setminus I}$; whereas for the
new basis the superscript varies with the index $v$ in the subscript.
More explicitly, the superscript is $\Delta\setminus J$ whenever $v\in
C^{J}$ where $W^{\Delta\setminus I}=\bigsqcup_{J\subseteq I} C^J$.
\erem

\begin{note}{\label{aproposds}} Henceforth throughout this paper we shall
fix the following notation: whenever $v\in C^{I}$ we shall denote
$f_v^{^{\Delta\setminus I}}$ simply by $f_{v}$. We can drop the
superscript in the notation without any ambiguity since $\{C^{I}:
I\subseteq \Delta\}$ are disjoint. Therefore with the modified
notation Lemma \ref{ds} implies that: $\{f_{v}: ~v\in
W^{\Delta\setminus I}=\bigsqcup_{J\subseteq I} C^J\}$ form an
$R(\tT)^{W}$-basis for $R(\tT)^{W_{\Delta\setminus I}}$ for every
$I\subseteq \Delta$ and
$$R(\tT)_{I}:=\bigoplus_{v\in C^I}R(\tT)^{W}\cdot f_v$$ satisfies
(1.13) and (1.14).
\end{note}

\subsubsection{Comparison with Topological K-theory}

Let $T_{comp}\subset T$ denote the maximal compact torus of $T$.
Then any complex algebraic $T$-variety can be viewed as a
topological $T_{comp}$-space. In particular, we have the algebraic
$K$-group $K_{T}(X)$ and the topological $K$-group
$K^{top}_{T_{comp}}(X)$. Now, since any algebraic vector bundle may
be regarded as a topological vector bundle we have a natural
homomorphism $K_{T}(X)\ra K^{top}_{T_{comp}}(X)$ (see pp.272
\cite{CG}).

\blem Let $X$ be a smooth projective variety on which $T$-acts with
finitely many fixed points. Then the canonical map $K_{T}(X)\ra
K^{top}_{T_{comp}}(X)$ is an isomorphism.\elem

{\it Proof:} The lemma follows by Proposition 5.5.6 of \cite{CG} since
$X\ra pt$ is a $T$-equivariant cellular fibration and
$K_{T}(pt)=R(T)\simeq R(T_{comp})=K^{top}_{T_{comp}}(pt)$.
\hfill $\Box$

\brem\label{top} Let $G_{comp}$ be a maximal compact subgroup of $G$
such that $T_{comp}=G_{comp}\cap T$ is a maximal torus in $G_{comp}$.
It has been proved in Theorem 4.4 of \cite{McL} that
$K^{top}_{G_{comp}}(X)\simeq (K^{top}_{T_{comp}}(X))^{W}$. Now, since
$K_{T}(X)\ra K_{T_{comp}}^{top}(X)$ is $W$-invariant, it further
follows from Theorem \ref{inv} that $K_{G}(X)\simeq
K^{top}_{G_{comp}}(X)$. \erem

\section{K-theory of regular embeddings}

We shall henceforth denote by $X$ a projective regular
compactification of $G$.  We follow the notations of \S1.1 together
with the following:

Let $\bar T$ denote the closure of $T$ in $X$. On $G$ the restriction
of the action of $diag(T)$ is given by $(t,t)\cdot g=tgt^{-1}$ for all
$g\in G$ and $t\in T$. This extends to an action on $X$. Thus $\bar T$
is an irreducible component of the fixed points of the torus $diag(T)$
and is therefore smooth (see Lemma 5.11.1. of \cite{CG}). Thus for the
left action of $T$ (i.e. for the action of $T\times \{1\}$), $\bar{T}$
is a smooth complete toric variety.

We now recall certain facts and notations from \S3.1 of \cite{Br}
suitably adapted to the setting of $K$-theory.

By Prop. A1 of \cite{Br}, $X^{T\times T}$ is contained in the union
$X_c$ of all closed $G\times G$-orbits in $X$; moreover all such
orbits are isomorphic to $G/B^{-}\times G/B$.  Therefore by Theorem 2
of \cite{VV}, $K_{T\times T}(X)$ embeds into $K_{T\times T}(X_c)$, the
latter being a product of copies of the ring $K_{T\times
T}(G/B^{-}\times G/B)$.

Let $\cf$ be the fan associated to $\bar T$ in $X_{*}(T)\otimes
\br$. Since $\bar{T}$ is complete, $\cf$ is a subdivision of
$X_{*}(T)\otimes \br$.  Moreover, since $\bar T$ is invariant under
$diag(W)$, the fan $\cf$ is invariant under $W$, too. Since $X$ is a
regular embedding, by Prop. A2 of \cite{Br}, it follows that
$\cf=W\cf_{+}$ where $\cf_{+}$ is the subdivision of the positive Weyl
chamber formed by the cones in $\cf$ contained in this
chamber. Therefore $\cf$ is a smooth subdivision of the fan associated
to the Weyl chambers, and the Weyl group $W$ acts on $\cf$ by
reflection about the Weyl chambers.  Let $\bar{T}^{+}$ denote the
toric variety associated to the fan ${\cal F}_{+}$.  Let ${\cf}(l)$
denote the set of maximal cones of $\cf$. Then we know that
${\cf}_{+}(l)$ parameterizes the closed $G\times G$-orbits in $X$.
Hence $X^{T\times T}$ is parametrized by ${\cf}_{+}(l)\times W\times
W$.  (The above facts follow from Prop. A1 and Prop. A2 of \cite{Br}.)

For $\sigma\in {\cf}_{+}(l)$, we denote by $Z_{\sigma}\simeq G/B^{-}\times
G/B$ the corresponding closed orbit with base point $z_{\sigma}$, and by

$$\iota_{\sigma}: K_{T\times T}(X) \ra K_{T\times
T}(Z_{\sigma})=K_{T\times T}(G/B^{-}\times G/B)$$

the restriction map. Moreover, for $f\in K_{T\times T}(Z_{\sigma})$
and $u,v\in W$, we denote by $f_{u,v}$, the restriction of $f$ to the
point $(u,v)z_{\sigma}$.

With the above notations, we have the following theorem.  For the
analogous result in the case of Chow ring see pp. 159 of \cite{Br}.

\bth \label{kring}
For any projective regular embedding $X$ of $G$, the map

$$\prod _{\sigma\in{\cf}_{+}(l)}\iota_{\sigma}: K_{T\times T}(X) \ra
\prod_{\sigma\in{\cf}_{+}(l)} K_{T\times
T}(G/B^{-}\times G/B)$$ is injective and its image consists in all
families $(f_{\sigma})$ $(\sigma\in {\cf}_{+}(l))$ in $R(T)\otimes R(T)$,
such that

\begin{enumerate}

\item[(i)] $f_{\sigma,us_{\alpha},vs_{\alpha}}\equiv f_{\sigma,
u,v}\pmod {(1-e^{-u(\alpha)}\otimes e^{-v(\alpha)})}$ whenever
$\alpha\in\Delta$ and the cone $\sigma\in{\cf}_{+}(l)$ has a facet
orthogonal to $\alpha$, and that

\item[(ii)] $f_{\sigma,u,v}\equiv f_{\sigma^{\prime}, u,v}\pmod
{(1-e^{-\chi})}$ whenever $\chi\in X^{*}(T)$ and the cones $\sigma$
and $\sigma^{\prime}\in{\cf}_{+}(l)$ have a common facet orthogonal to
$\chi$.

\end{enumerate}
(In $(ii)$, $\chi$ is viewed as a character of $T\times T$ which is
trivial on $diag(T)$ and hence is a character of $T$.)

\eeth

{\it Proof:} In the proof of the Theorem on pp.160 of \cite{Br} we
have a complete description of all $T\times T$-invariant irreducible
curves in $X$. We briefly recall here this description.

{\it $T\times T$-invariant curves in $X$:} Let $\gamma$ be
a $T\times T$-invariant irreducible curve in $X$. Then $\gamma$ joins
two $T\times T$-fixed points in $X$ and one of the following cases
occur:

{\bf (1)} $\gamma$ lies inside a closed orbit $Z_{\sigma}$. Thus by
    the description of $T$-invariant curves in $G/B$ (see \S6.5 of
    \cite{Br2}) it follows that $\gamma$ is conjugate in $W\times W$
    to a curve $\gamma^{\prime}$ joining $z_{\sigma}$ to
    $(s_{\alpha},1)z_{\sigma}$ or to $(1,s_{\alpha})z_{\sigma}$ where
    $z_{\sigma}$ is the base point of $Z_{\sigma}$.

{\bf (2)} $\gamma$ is conjugate in $W\times W$ to a curve $\gamma^{\prime}$
    joining the $T\times T$-fixed points $z_{\sigma}$ and
    $(s_{\alpha},s_{\alpha})z_{\sigma}$ of the closed orbit
    $Z_{\sigma}$, where $\gamma^{\prime}$ is not contained in
    $Z_{\sigma}$. In this case the cone $\sigma$ in $\cf_{+}(l)$ has a
    facet orthogonal to $\alpha$.

{\bf (3)} $\gamma$ is conjugate in $W\times W$ to a projective line
    $\gamma^{\prime}$ joining the $T\times T$-fixed points
    $z_{\sigma}$ and $z_{\sigma^{\prime}}$ which are respectively the
    base points of distinct closed orbits $Z_{\sigma}$ and
    $Z_{\sigma^{\prime}}$. In this case the cones $\sigma$ and
    $\sigma^{\prime}$ in $\cf_{+}(l)$ have a common facet.

In particular we observe that the set of $T\times T$-invariant
irreducible curves in $X$ is finite.

Therefore by Theorem \ref{im}, the image of

$$\iota^{*}: K_{T\times T}(X) \ra K_{T\times T}(X^{T\times T})$$

\noindent
is defined by linear congruences $f_x \equiv f_y \pmod
{(1-e^{-\chi})}$ whenever $x,y\in X^{T\times T}$ are connected by a
curve where $T\times T$ acts by the character $\chi$.

Further, observe that $T\times T$ acts on the curve joining
$z_{\sigma}$ to $(s_{\alpha},s_{\alpha})z_{\sigma}$ by the character
$(\alpha,\alpha)$, and on the curve joining $z_{\sigma}$ to
$z_{\sigma^{\prime}}$ by the character $\chi$ where $\sigma$ and
$\sigma^{\prime}$ have a common facet orthogonal to $\chi$.  It
therefore follows that the curves of type (1) define the image of
$\prod _{\sigma\in{\cf}_{+}}\iota_{\sigma}$, whereas curves of type
(2) and (3) lead to congruences (i) and (ii). \hfill $\Box$

\bcor{\label{co1}} The ring $K_{G\times G}(X)$ consists in all families
$(f_{\sigma}) (\sigma\in \cf_{+}(l))$ of elements of $R(T)\otimes
R(T)$ such that:

\begin{enumerate}
\item[(i)] $(s_{\alpha},s_{\alpha})f_{\sigma}\equiv f_{\sigma} \pmod
{(1-e^{-\alpha}\otimes e^{-\alpha})}$ whenever $\alpha\in\Delta$ and
the cone $\sigma\in{\cf}_{+}(l)$ has a facet orthogonal to $\alpha$,
and that

\item[(ii)] $f_{\sigma}\equiv f_{\sigma^{\prime}}\pmod
{(1-e^{-\chi})}$ whenever $\chi\in X^{*}(T)$ and the cones $\sigma$ and
$\sigma^{\prime}\in{\cf}_{+}(l)$ have a common facet orthogonal to $\chi$.

\end{enumerate}
\ecor {\it Proof:} By the isomorphism $(b)$ in \S1.1, the ring
$K_{G\times G}(G/B^{-}\times G/B)$ is isomorphic to $K_{T\times
T}(G/B^{-}\times G/B))^{W\times W}$. It is further isomorphic to
$R(T)\otimes R(T)$ via restriction to $z_{\sigma}$. Moreover,
restriction of $f\in (K_{T\times T}(G/B^{-}\times G/B))^{W\times
W}\simeq R(T)\otimes R(T)$ to $(u,v)z_{\sigma}$ is equal to
$(u,v)f_{\sigma}$ where $f_{\sigma}$ denotes the restriction of $f$ to
$z_{\sigma}$. So the relations (i) and (ii) of Theorem \ref{kring}
reduce to (i) and (ii) of Cor.\ref{co1}. \hfill $\Box$

We have the following relation between $K_{G\times G}(X)$ and
$K_{T\times T}(\bar{T})$. This is analogous to the relation for
semisimple adjoint groups and equivariant cohomology, due to
Littelmann and Procesi (see \cite{LP}), and to the corresponding
relation for the equivariant Chow ring of a regular group
compactification due to Brion (see Cor.2 in \S3.1 of \cite{Br}).

\bcor{\label{co2}} The inclusion $\bar{T}\hra X$ induces the following
isomorphisms:
$$K_{G\times G}(X)\simeq K_{T\times T}(\bar T)^{W}\simeq
(K_{T}(\bar{T})\otimes R(T))^W$$ where the $W$-action on $K_{T\times
T}(\bar T)$ is induced from the action of $diag(W)$ on $\bar{T}$.

\ecor

{\it Proof:} Let $N$ be the normalizer of $T$ in $G$ and let $\bar{N}$
be its closure in $X$. Observe that $\bar N$ is the disjoint union of
$(w,1)\bar T$ for $w\in W$. This can be seen as follows:
We have $N=\bigcup_{w\in W} w T.$ This implies that
$$\bar{N}=\bigcup_{w\in W}(w,1) \bar{T}.$$

Further, the map $y\mapsto (w,1)y$ $~\forall ~y\in \bar{T}$ is an
isomorphism from $\bar{T}$ to $(w,1)\bar{T}$ on which the $T\times
T$-action is twisted by $(w,1)$. In particular, the $T\times T$-fixed
points in $(w,1)\bar{T}$ are $(w,1)\cdot\bar{T}^{T\times T}$.

Now, the set of fixed points $\bar{T}^{T\times T}$ is parametrized by
$\cf_{+}(l)\times diag(W)$. Therefore the set of $T\times T$-fixed
points $(w,1)\cdot\bar{T}^{T\times T}$ is parametrized by
$\cf_{+}(l)\times (w,1) diag(W)$.  However we know that $X^{T\times
T}$ is parametrized by $\cf_{+}(l)\times W\times W$ where $W\times
W=\bigcup_{w\in W}(w,1)diag(W)$.

It follows that $(w,1)\bar{T}$ are disjoint, for otherwise the
intersection of two of them should contain $T\times T$-fixed points
which is a contradiction. Therefore we have:

$$\bar{N}=\bigsqcup_{w\in W}(w,1)\bar{T}\eqno(2.1)$$

where for each $w\in W$, $(w,1)\bar{T}$ is an irreducible variety
isomorphic to $\bar{T}$ with the appropriate twist for the $T\times
T$-action.

In particular, $\bar N$ contains all fixed points of $T\times T$. It
follows that restriction
$$K_{T\times T}(X)\ra K_{T\times T}(\bar{N})$$ is injective.

Further, taking invariants of $K_{T\times T}(X)$ under $W\times
W$, we see that the relations arising from curves of type (1) reduce
to $(s_{\alpha},1)(f_{\sigma})\equiv f_{\sigma}\pmod
{(1-e^{-\alpha})\otimes 1}$ or $(1,s_{\alpha})(f_{\sigma})\equiv
f_{\sigma}\pmod {1\otimes (1-e^{-\alpha})}$ for $f_{\sigma}\in
R(T)\otimes R(T)$ for every $\sigma\in {\cf}_{+}(l)$.

However these relations trivially hold in $R(T)\otimes R(T)$, due
to the fact that $s_{\alpha}(e^{\lambda})-e^{\lambda}$ is divisible in
$R(T)$ by the element $1-e^{-\alpha}$ for every $\alpha\in \Delta$
and $\lambda\in \wt{\Lambda}$.

Therefore the non-trivial relations which describe the image of
$K_{T\times T}(X)^{W\times W}$ arise from the curves of type (2) and
(3).

From the description of $T\times T$-invariant curves in $X$ of type
(2) and (3), it follows that any curve of type $(2)$ or $(3)$ lies in
$(w,1)\bar{T}$ for a unique $w\in W$ (since any such curve is
conjugate in $W\times W$ to a curve lying in $\bar{T}$). Thus $\bar N$
contains all $T\times T$-invariant curves which are not in any closed
$(G\times G)$-orbit, that is curves of type (2) and (3). Thus we see
that the restriction to $\bar N$ induces an isomorphism

$$K_{T\times T}(X)^{W\times W}\simeq K_{T\times T}(\bar N)^{W\times
W}.$$

Further, by $(2.1)$ it follows that

$$ K_{T\times T}(\bar{N})\simeq \bigoplus_{w\in W}(w,1)
K_{T\times T}(\bar{T})$$ where $(w,1)$ denotes the isomorphism on
$K$-rings induced by the above isomorphism from $\bar{T}$ to
$(w,1)\bar{T}$.

Thus the $W\times W$-module structure on $K_{T\times T}(\bar{N})$ is
induced from the $diag(W)$-module structure on $K_{T\times
T}(\bar{T})$. Thus we have
$$K_{T\times T}(\bar N)^{W\times W}\simeq K_{T\times T}(\bar
T)^{W}.$$

Therefore we have the following isomorphisms

$$K_{T\times T}(X)^{W\times W}\simeq K_{T\times T}(\bar N)^{W\times
W}\simeq K_{T\times T}(\bar T)^{W}\simeq (K_{T}(\bar{T})\otimes
R(T))^W.$$

\noindent
(The last isomorphism is a consequence of the fact that we
have a split exact sequence
$$1\ra diag(T)\ra T\times T\ra T\ra 1$$ where the second map
is $(t_1,t_2)\ra t_1 t_2^{-1}$, and the splitting is given
by $t\ra (t,1)$.

Thus $T\times T$ is canonically isomorphic to $diag(T)\times
(T\times \{1\})$. Furthermore, by the definition of $T\times T$
action on $\bar T$ we see that $diag(T)$ acts trivially on $\bar
T$. Therefore we have a ring isomorphism $K_{T\times T}(\bar
T)\simeq R(diag(T)) \otimes K_{T}(\bar T)$ (see 5.2.4 pp.244 of
\cite{CG}). This isomorphism is further $W$-invariant since the
$W$-action on the $K$-rings is induced from the action of $diag(W)$ on
$\bar T$.) \hfill $\Box$

\brem\label{prop} Since $\cf_{+}$ is a subdivision of the fan
associated to the positive Weyl chamber, we have an induced proper
morphism $\bar{T}^{+}\ra \ba^l$. (Here $T$ acts linearly on $\ba^l$
with weights being the simple roots.) Therefore, by the valuative
criterion for properness it follows that $T$ acts on $\bar{T}^{+}$
with {\it enough limits} (see pp. 19 of \cite{VV} for the definition
of action with enough limits). Thus by Cor. 5.11 and Cor. 5.12 of
\cite{VV} it further follows that the restriction homomorphism
$$K_{T}(\bar{T}^+)\ra \prod_{\sigma\in\cf_+(l)}
R(T_{\sigma})=K_T((\bar{T}^{+})^{T})$$ is injective and an element
$(a_{\sigma})\in\prod_{\sigma}R(T_{\sigma})$ is in the image of this
homomorphism if and only if for any two {\it adjacent} maximal cones
$\sigma$ and $\sigma^{\prime}$, the restrictions of $a_{\sigma}$ and
$a_{\sigma^{\prime}}$ to $R(T_{\sigma\cap\sigma^{\prime}})$ coincide,
where $T_{\tau}\subseteq T$ denotes the stabilizer along the orbit
$O_{\tau}$ for every $\tau\in \cf_{+}$. Further, since $\bar{T}^{+}$
is smooth, and $T$ acts on $\bar{T}^{+}$ with finitely many fixed
points and finitely many invariant curves, it can be seen (see
Cor 5.11 of \cite{VV}) that Theorem \ref{im} holds for
$\bar{T}^{+}$. More precisely, the image of $K_{T}(\bar{T}^{+})$
consists of elements $(a_{\sigma})\in \prod_{\sigma}R(T_{\sigma})$
such that $a_{\sigma}-a_{\sigma^{\prime}}\equiv 0 \pmod
{(1-e^{-\chi})}$ whenever $\sigma$ and $\sigma^{\prime}$ have a common
facet orthogonal to $\chi\in X^{*}(T)$.  \erem

The following proposition is a consequence of Corollary. \ref{co1}.

\bpropo\label{key} We have the following chain of inclusions as
$R(G)\otimes R(G)$-modules:
$$R(T)\otimes R(G)\subseteq K_{T}(\bar{T}^{+})\otimes R(G)\subseteq
K_{G\times G}(X) \subseteq R(T)^{|\cf_{+}(l)|}\otimes R(T).$$

Moreover, $K_{G\times G}(X)$ is a module over $R(T)\otimes
R(G)$. \epropo

{\it Proof:} From the split exact sequence $$1\ra diag(T)\ra
T\times T\ra T\ra 1$$ it follows that $R(T)\otimes
R(T)\simeq R(T\times \{1\})\otimes R(diag(T))$.

Recall that
$$\bar{N}=\bigsqcup_{w\in W}(w,1) \bar{T},$$ and any $T\times
T$-invariant curve of type (2) or (3) in $X$ lies in $(w,1)\bar{T}$
for some $w\in W$.  In particular, it follows that $diag(T)$ acts
trivially on the curve $\gamma$ joining $(w,1)z_{\sigma}$ and
$(w,1)(s_{\alpha},s_{\alpha})z_{\sigma}$ for every
$\sigma\in\cf_{+}(l)$ having a facet orthogonal to $\alpha\in\Delta$,
and $w\in W$. Moreover, $T\times \{1\}$ acts on $\gamma$ by the
character $\alpha$ (where the action is twisted by $(w,1)$ on the
curves lying in $(w,1)\bar{T}$).

Similarly, $diag(T)$ acts trivially on the curve $\gamma$ joining
$(w,1)z_{\sigma}$ and $(w,1)z_{\sigma^{\prime}}$, and $T\times \{1\}$
acts on $\gamma$ by the character $\chi$, for all cones $\sigma$ and
$\sigma^{\prime}\in \cf_{+}(l)$ having a common facet orthogonal to
$\chi$.

Hence by Cor.\ref{co1} it follows that $K_{G\times G}(X)$ consists
in all families $(f_{\sigma})(\sigma\in\cf_{+}(l))$ of elements of $
R({T\times\{1\}})\otimes R(diag(T)$ such that:

\begin{enumerate}
\item[$(i)$] $(1,s_\alpha)f_{\sigma}(u,v)\equiv f_{\sigma}(u,v)\pmod
{(1-e^{-\alpha(u)})}$ whenever $\alpha\in\Delta$ and the cone
$\sigma\in\cf_{+}(l)$ has a facet orthogonal to $\alpha$.

\item[$(ii)$] $f_{\sigma}\equiv f_{\sigma^{\prime}}\pmod
{(1-e^{-\chi(u)})}$ whenever $\chi\in X^{*}(T)$ and the cones $\sigma$
and $\sigma^{\prime}\in{\cf}_{+}(l)$ have a common facet orthogonal to
$\chi$.
\end{enumerate}
where $u$ and $v$ denote the variables corresponding to
$R(T\times\{1\})$ and $R(diag(T))$ respectively.


Now, by Remark \ref{prop} it follows that $K_{T}(\bar{T}^{+})\otimes
R(T)^{W}\subseteq \prod_{\sigma\in\cf_{+}(l)} R(T_{\sigma})\otimes
R(T)$ is generated by the elements $(a_{\sigma})\otimes b$, where
$a_{\sigma}-a_{\sigma^{\prime}}\equiv 0 \pmod{(1- e^{-\chi})}$,
whenever $\sigma$ and $\sigma^{\prime}$ share a facet orthogonal to
$\chi\in X^*(T)$, and $b\in R(T)^{W}$.

Therefore, by identifying both $T\times\{1\}$ and $diag(T)$ naturally
with $T$ keeping track of the ordering, we see that
$K_{T}(\bar{T}^{+})\otimes R({T})^{W}$ satisfies the relations $(i)$
and $(ii)$.  Therefore it is a submodule of $K_{G\times
G}(X)$. Moreover, since $K_{T}(\bar{T}^{+})$ is an algebra over
$R(T)$, it follows that $K_{G\times G}(X)$ is a module over
$R(T)\otimes R(G)$.  \hfill $\Box$

We give below the example of wonderful compactification of
$PGL(2,\bc)$. In particular we shall clearly see the curves of type
(1) and (2) in this case.

\beg Let $G=PGL(2,\bc)=SL(2,\bc)/{\pm Id}$. Then the projective space
$\bp(M(2,\bc))$ is the wonderful compactification of $PGL(2,\bc)$, on
which the action of $PGL(2,\bc)\times PGL(2,\bc)$ by
multiplication on the left and on the right extends.

Let $E_{ij}$ denote the elementary matrix with $1$ as $(i,j)th$ entry
and $0$ elsewhere for $1\leq i,j\leq 2$. In this case the Weyl group
is $W=\{1=Id,s_{\alpha}=-E_{12}+E_{21}\}$, and $\bar{T}\simeq \bp^1$
consists of the diagonal matrices in $\bp(M(2,\bc))$.

Further, the unique closed $PGL(2,\bc)\times PGL(2,\bc)$-orbit
consists of the matrices of rank $1$ in $ \bp(M(2,\bc))$ and is
isomorphic to $PGL(2,\bc)\times PGL(2,\bc)/(B^{-}\times B^{+})$,
choosing as base point the matrix $E_{11}$. Furthermore, $PGL(2,\bc)$
is the open orbit with base point $Id$.

The four $T\times T$ fixed points of $\bp(M(2,\bc))$ are: $E_{11}$,
$E_{12}=(1,s_{\alpha})E_{11}$, $E_{21}=(s_{\alpha},1)E_{11}$ and
$E_{22}=(s_{\alpha},s_{\alpha})E_{11}$.  Further, the  $T\times T$
curves are the following:

\begin{enumerate}
\item[(1)] $aE_{11}+bE_{12}$; $aE_{11}+bE_{21}$;
$aE_{12}+bE_{22}$; $aE_{21}+bE_{22}$.

\item[(2)] $aE_{11}+bE_{22}$ and $aE_{12}+bE_{21}$.

\end{enumerate}
where $aE_{ij}+bE_{pq} ~~\forall~a,b\in\bc$, denotes the projective
line joining $E_{ij}$ and $E_{pq}~$ in $\bp(M(2,\bc))$ for $i,j,p,q\in
\{1, 2\}$. Pictorially we can view these curves as follows:

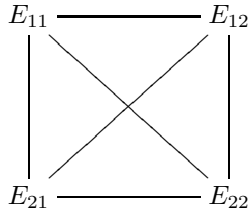
\[
\xymatrix{
E_{11} \ar@{-}[rr] & & \ar@{-}[dd] \ar@{-}[ddll] E_{12} \\
& & \\
E_{21} \ar@{-}[uu] & & \ar@{-}[ll] \ar@{-} [uull] E_{22} \\
}
\]

Thus we see that the curves of type (1) lie entirely in the unique
closed orbit, whereas the curves of type (2) meet the open orbit.

Moreover, $\bar{N}=\bar{T}\sqcup (s_{\alpha},1)\bar{T}$ is the union of
diagonal and the antidiagonal matrices. Hence $\bar{N}$ contains
only the curves of type (2) and does not contain the curves of type
(1).

In this case we do not have curves of type (3) since there is a unique
closed $G\times G$-orbit.  \eeg

\brem\label{ext1} Note that {\it all the results in this section hold
analogously for $K_{\tG\times \tG}(X)$ and $K_{\tT\times \tT}(X)$}
where we take the natural actions of $\tG\times \tG$ and $\tT\times
\tT$ through the canonical surjections to $G\times G$ and $T\times T$
respectively.  \erem

\subsection{Determination of the structure of $K_{G\times G}(X)$}

Let $X:=\bar G$ be a projective regular embedding of $G$ and let $\bar
T$ be the corresponding torus embedding.

Let $\cf$ be the (smooth projective) fan associated to $\bar
T$. Recall that the Weyl group $W$ acts on $\cf$ by reflection about
the Weyl chambers and the cones in $\cf$ get permuted by this action
of $W$, and each cone is stabilized by the reflections corresponding
to the walls of the Weyl chambers on which it lies. Let $W_{\tau}$
denote the subgroup of $W$ which fixes the cone $\tau\in\cf$. Then in
particular, $W_{\sigma}=\{1\}~\forall~ \sigma\in\cf(l)$, and
$W_{\{0\}}=W$.

Let $\{\rho_j:j=1,\ldots d\}$ denote the set of edges of the fan $\cf$
and let $\tau(1)$ denote the set of edges of the cone $\tau$ for every
$\tau\in\cf$. Let $v_j$ denote the primitive vector along the edge
$\rho_j$. Let $O_{\tau}$ denote the $T$-orbit in $\bar T$
corresponding to $\tau\in\cf$. Let $L_{j}$ denote the $T$-equivariant
line bundle on $\bar{T}$ corresponding to the edge $\rho_{j}$. We note
that, $L_{j}$ has a $T$-invariant section $s_j$ whose zero locus is
$\bar{O_{{\rho}_j}}$.  Recall that $\bar{T}^{+}$ denotes the toric
variety associated to the fan $\cf_{+}$.

Let $X_F:=\prod_{\rho_j\in F}(1- X_j)$ in the Laurent polynomial
algebra $\bz[X_1^{\pm 1},\ldots X_d^{\pm 1}]$, for every $F\subseteq
\{\rho_j:j=1,\ldots d\}$. In particular, let
$X_{\tau}:=X_{\tau(1)}=\prod_{\rho_j\in \tau(1)} (1-X_j)$ for every
$\tau\in\cf$.

Recall from Theorem 6.4 of \cite{VV} we have the following
Stanley-Reisner presentation of the $T$-equivariant $K$-ring of
$\bar{T^{+}}$:
$$K_{T}(\bar T^{+})\simeq\bz[X_j^{\pm 1}: \rho_j\in\cf_{+}(1)]/\langle
X_F, ~~\forall~ F\notin \cf_{+} \rangle $$ where under the above
isomorphism $X_j$ maps to $[L_j]$.

Further, we have the additive decomposition
$K_{T}(\bar T^{+})=\bigoplus_{\tau\in \cf_{+}} C_{\tau}$, where
$$C_{\tau}:=X_{\tau}\cdot \bz[X_j^{\pm 1}:\rho_j \in \tau(1)].$$ Since
we do not have an immediate reference for the above additive
decomposition which may be well known, we give a proof of it in the
following lemma.

\blem{\label{addec}} We have the additive decomposition $K_{T}(\bar
T^{+})=\bigoplus_{\tau\in \cf_{+}} C_{\tau}$, where
$$C_{\tau}:=X_{\tau}\cdot \bz[X_j^{\pm 1}:\rho_j \in \tau(1)].$$ \elem

{\it Proof:} Let $\ci$ be any finite indexing set. Then $\bz[X_j^{\pm
1}: j\in \ci]=R[X_i^{\pm 1}]=R\bigoplus (1-X_i)R[X_i^{\pm 1}]$ where
$R:=\bz[X_j^{\pm 1}: j\in \ci, j\neq i]$.  By induction on $|\ci|$ we
have the following direct sum decomposition: $$\bz[X_j^{\pm 1}:
j\in\ci]=\bigoplus_{F\subseteq \ci} X_F\cdot \bz[X_j^{\pm 1}: j\in
F]$$ where $X_{F}:=\prod_{j\in F}(1-X_j)$.

Now, let $F\subseteq \ci$ and $i\notin F$. Then we have $X_F\cdot
\bz[X_j^{\pm 1}: j\in\ci]=X_{F}\cdot R[X_i^{\pm 1}]=X_{F}\cdot
R\bigoplus X_{F^{\prime}}\cdot R[X_i^{\pm 1}]$ where $R:=\bz[X_j^{\pm
1}: j\in\ci, j\neq i]$ and $X_{F^{\prime}}=X_F(1-X_i)$.  Thus by
induction on $|\{i:i\notin F\}|$ it follows that we have the following
direct sum decomposition: $$X_F\cdot \bz[X_j^{\pm 1}: j\in\ci]=
\bigoplus_{F\subseteq F^{\prime}} X_{F^{\prime}}\cdot \bz[X_j^{\pm 1}:
j\in F^{\prime}].$$

In particular, applying the above arguments for $\ci=\cf_{+}(1)$ we
get:
$$\bz[X_j^{\pm 1}: \rho_j\in\cf_{+}(1)]=\bigoplus_{F\subseteq
\cf_+(1)} X_F\cdot \bz[X_j^{\pm 1}: \rho_j\in F]$$ and further for
$F\subseteq \cf_{+}(1)$ such that $F\notin \cf_{+}$ we get:
$$X_F\cdot \bz[X_j^{\pm 1}: \rho_j\in\cf_{+}(1)]=
\bigoplus_{F\subseteq F^{\prime}} X_{F^{\prime}}\cdot \bz[X_j^{\pm 1}:
\rho_j\in F^{\prime}]$$

Since $F\subseteq F^{\prime}$ implies $F^{\prime}\notin \cf_{+}$ it
follows that we have the following direct sum decomposition:
$$\langle X_F, ~~\forall~ F\notin \cf_{+} \rangle= \bigoplus_{F\notin
\cf_+} X_F\cdot \bz[X_j^{\pm 1}: \rho_j\in F]$$

The lemma now follows by the Stanley-Reisner presentation of
$K_{T}(\bar T^{+})$.\hfill $\Box$

\brem{\label{toraddec}} Note that although we state Lemma \ref{addec}
for $\bar{T}^{+}$ it is not hard to see that an analogous additive
decomposition holds for the $T$-equivariant $K$-ring of any smooth
$T$-toric variety.  \erem

Whenever $\sigma,\tau\in \cf$ (resp. $\cf_{+}$) span a cone in $\cf$
(resp. $\cf_{+}$), we shall denote the cone spanned by them as
$\gamma:=\langle \tau,\sigma \rangle$.

\bth{\label{kdec}} Let $X:=\bar G$ be a projective regular embedding
of $G$ and let $\bar T$ be the corresponding torus embedding. Then,
$K_{G\times G}(X)$ has the following direct sum decomposition as
$1\otimes R(G)$-module:

$$K_{G\times G}(X) \simeq \bigoplus_{\tau\in {\cal F}_{+}}
C_{\tau}\otimes R(T)^{W_{\tau}} $$

where $R(G)=R(T)^{W}$ acts naturally on the second factor in each
piece of the above decomposition.  Further, the {\em multiplicative
structure} of $K_{G\times G}(X)$ can be described from the above
decomposition as follows: Let $a_{\tau}\otimes b_{\tau}\in
C_{\tau}\otimes R(T)^{W_{\tau}}$ and $a_{\sigma}\otimes b_{\sigma}\in
C_{\sigma}\otimes R(T)^{W_{\sigma}}$. Then
\[(a_{\tau}\otimes b_{\tau})\cdot(a_{\sigma}\otimes
b_{\sigma})=\left\{ \begin{array}{ll}a_{\tau}\cdot a_{\sigma}\otimes
b_{\tau}\cdot b_{\sigma}, &\mbox {if $\tau$ and $\sigma$ span the cone
$\gamma$}\\ 0 & \mbox {if $\tau$ and $\sigma$ do not span a cone in
$\cf_{+}$}\end{array}\right.\]

{\em(Note that $a_{\tau}\cdot a_{\sigma}\otimes b_{\tau}\cdot
b_{\sigma}\in C_{\gamma}\otimes R(T)^{W_{\gamma}}$, and the
multiplication in the first factor is as in $K_{T}(\bar T^{+})$ where
$C_{\tau} \cdot C_{\sigma}\subseteq C_{\gamma}$.)}
\eeth

\noindent
{\it Proof:} We have the following isomorphisms by Cor.\ref{co2}:

$$K_{G\times G}(X)\simeq K_{T\times T}(\bar T)^{W}\simeq
(K_{T}(\bar{T})\otimes R(T))^W$$

Now by Theorem 6.4 of \cite{VV} we have the following Stanley-Reisner
presentation of the $T$-equivariant $K$-ring of $\bar{T}$:

$$K_{T}(\bar T)\simeq\bz[X_1^{\pm 1},\ldots X_d^{\pm 1}]/\langle X_F,
~~\forall~ F\notin \cf \rangle \eqno(2.2)$$

Since $W$ acts on $\cf$, we have an action of $W$ on $\bz[X_1^{\pm
1},\ldots X_d^{\pm 1}]$, given by $w(X^{\pm 1}_{\rho_j})=X^{\pm
1}_{w(\rho_j)}$ for every $w\in W$. Therefore, $w(X_F)=X_{w(F)}$ for
$F\subseteq \{\rho_j:j=1,\ldots d\}$ and $w\in W$, and since $W$
permutes the cones of $\cf$ we further get an action of $W$ on the
Stanley-Reisner algebra $\bz[X_1^{\pm 1},\ldots X_d^{\pm 1}]/\langle
X_F ~\forall~ F\notin \cf \rangle$. The above isomorphism is an
isomorphism of $W$-modules, where the $W$-action on $K_{T}(\bar{T})$
is induced by the $diag(W)$-action on $\bar{T}$.

Further, (see Lemma \ref{addec} and Remark \ref{toraddec})
$$\bz[X_1^{\pm 1},\ldots X_d^{\pm 1}]/\langle X_F, ~~\forall~ F\notin
\cf \rangle=\bigoplus_{\tau\in\cf} X_{\tau}\cdot \bz[X_j^{\pm
1}:\rho_j \in \tau(1)]$$  where we have the
natural action of $W$ on the right hand side given by:
$$w\cdot (X_{\tau}\cdot \bz[X_{j}^{\pm 1}:\rho_{j}\in\tau(1)])=
X_{w(\tau)}\cdot \bz[X_{j}^{\pm 1}:\rho_{j}\in w(\tau)(1)],
~~~\forall~~w\in W.$$ Therefore we have:
$$K_{T}(\bar T)=\bigoplus_{\tau\in{\cal F}_{+}} \bigoplus_{w\in
W/W_{\tau}} X_{w(\tau)}\cdot \bz[w(X_j^{\pm 1}):\rho_j \in \tau(1)]$$

Here $W_{\tau}$ denotes the subgroup of $W$ which fixes the cone
$\tau\in{\cal F}_{+}$. Hence we have as $W$-modules:

$$K_{T}(\bar T)= \bigoplus_{\tau\in{\cal F}_{+}}
Ind^{W}_{W_{\tau}}C_{\tau}$$ where, $C_{\tau}:=X_{\tau}\cdot
\bz[X_j^{\pm 1}:\rho_j \in \tau(1)]$.

Further, since $C_{\tau}$ is fixed by $W_{\tau}$, hence
$Ind^{W}_{W_{\tau}}C_{\tau}\simeq \bz[W/W_{\tau}]\otimes
C_{\tau}$.

Thus we have:
$$K_{T}(\bar T)\otimes R(T)= \bigoplus_{\tau\in {\cal F}_{+}}
\bz[W/W_{\tau}]\otimes C_{\tau}\otimes R(T).$$

Now by taking $W$-invariants on either side we get:

$$(K_{T}(\bar{T})\otimes R(T))^W= \bigoplus_{\tau\in {\cal F}_{+}}
C_{\tau}\otimes R(T)^{W_{\tau}}. \eqno(2.3)$$

Thus we get the following additive decomposition:
$$K_{G\times G}(X) \simeq \bigoplus_{\tau\in {\cal F}_{+}}
C_{\tau}\otimes R(T)^{W_{\tau}}. $$

We shall now describe the multiplication on the right hand side of the
above decomposition which will make the above isomorphism a ring
isomorphism.

First we note that the isomorphism $K_{G\times G}(X)\simeq
(K_{T}(\bar{T})\otimes R(T))^W $ is a ring isomorphism and hence
preserves the multiplicative structure.

Further, since $(2.2)$ is a ring isomorphism the multiplication in
$K_{T}(\bar T)$ is determined by multiplication in the Stanley-Reisner
ring $\bz[X_1^{\pm 1},\ldots X_d^{\pm 1}]/\langle X_F, ~~\forall~
F\notin \cf \rangle$ (see \S6.2 of \cite{VV}).  Moreover, the
multiplication in
$$\bz[X_1^{\pm 1},\ldots X_d^{\pm 1}]/\langle X_F, ~~\forall~ F\notin
\cf \rangle =\bigoplus_{\tau\in\cf}C_{\tau}$$ is determined by
the products $C_{\tau}\cdot C_{\sigma}$ for
$\tau,\sigma\in\cf$, where $C_{\tau}\cdot C_{\sigma}\subseteq
C_{\gamma}$ for $\gamma=\langle \tau, \sigma \rangle$.  Similarly, the
multiplication in $K_{T}(\bar
T^{+})=\bigoplus_{\tau\in\cf_{+}}C_{\tau}$ is defined by multiplying
$C_{\tau}$ and $C_{\sigma}$ for $\tau$ and $\sigma$ in
$\cf_+$. Moreover, if $\tau$ and $\sigma$ span a cone in $\cf_{+}$
then $C_{\tau}\cdot C_{\sigma}\subseteq C_{\gamma}$ where
$\gamma=\langle \sigma, \tau \rangle$, and if $\tau$ and $\sigma$ do
not span any cone in $\cf_{+}$ then $C_{\tau}\cdot C_{\sigma}=0$.

Furthermore, the multiplicative structure on $K_{T}(\bar T)=
\bigoplus_{\tau\in{\cal F}_{+}} Ind^{W}_{W_{\tau}}C_{\tau}$ is induced
from the multiplication in $K_{T}(\bar T^{+})=\bigoplus_{\tau\in\cf_+}
C_{\tau}$. This is because, for $w^{\prime},w^{\prime\prime}\in W$,
$w^{\prime}(\sigma)$ and $w^{\prime\prime}(\tau)$ span a cone in $\cf$
if and only if $\tau$ and $\sigma$ span a cone in $\cf_{+}$, and there
exists $w\in W$ such that $w(\sigma)=w^{\prime}(\sigma)$ and
$w(\tau)=w^{\prime\prime}(\tau)$. Moreover, if $\tau$ and $\sigma$
span $\gamma$ in $\cf_+$, $w(\tau)$ and $w(\sigma)$ span $w(\gamma)$
in $w(\cf_+)$ for every $w\in W$.

In particular, let $g_{\sigma}=w^{\prime}(f_{\sigma})\in
C_{w^{\prime}(\sigma)}$ and $g_{\tau}=w^{\prime\prime}(f_{\tau})\in
C_{w^{\prime\prime}(\tau)}$, where $f_{\sigma}\in C_{\sigma}$ and
$f_{\tau}\in C_{\tau}$ for $\tau, \sigma \in \cf_+$. Then
$g_{\sigma}\cdot g_{\tau}=0$ if $\sigma$ and $\tau$ do not span a cone
in $\cf_{+}$, or if there does not exist any $w\in W$ such that
$w(\sigma)=w^{\prime}(\sigma)$ and $w(\tau)=w^{\prime\prime}(\tau)$.
Otherwise $g_{\sigma}\cdot g_{\tau}=w(f_{\sigma}\cdot f_{\tau})\in
C_{w(\gamma)}$, where $\gamma=\langle \tau, \sigma \rangle$, and
$w(\sigma)=w^{\prime}(\sigma)$ and $w(\tau)=w^{\prime\prime}(\tau)$
for $w\in W$.

Further note that, whenever $\gamma=\langle \tau, \sigma \rangle$ in
$\cf_{+}$ we have the product: $R(T)^{W_{\tau}}\cdot
R(T)^{W_{\sigma}}\subseteq R(T)^{W_{\gamma}}$, where $R(T)^{W_{\tau}}$
and $R(T)^{W_{\sigma}}$ are both subrings of $R(T)^{W_{\gamma}}$.

Thus we see that the identity $(2.3)$ above induces a multiplicative
isomorphism, where the multiplication in $\bigoplus_{\tau\in {\cal
F}_{+}} C_{\tau}\otimes R(T)^{W_{\tau}}$ is as described in the
statement of the theorem. \hfill $\Box$

\bcor\label{kfilt}
The ring $K_{G\times G}(X) \simeq \bigoplus_{\tau\in {\cal F}_{+}}
C_{\tau}\otimes R(T)^{W_{\tau}} $ admits a multifiltration
$\{F_{\tau}\}_{\tau\in\cf_{+}}$, where the filtered pieces are
$$F_{\tau}=\bigoplus_{\tau\prec\sigma} C_{\sigma}\otimes
R(T)^{W_{\sigma}},$$ where $F_{\tau}\supseteq F_{\sigma}$ whenever
$\tau\prec\sigma$, and $F_{\{0\}}= K_{G\times G}(X)$.  Further, under
the multiplication described in Theorem \ref{kdec}, we have
$F_{\tau}\cdot F_{\sigma}\subseteq F_{\gamma}$ where $\gamma=\langle
\tau,\sigma\rangle$. In particular, $F_{\{0\}}\cdot F_{\tau}\subseteq
F_{\tau}$ for all $\tau\in\cf_{+}$. Moreover, since
$K_{T}(\bar{T}^{+})\otimes R(T)^{W} \subseteq F_{\{0\}}$, it follows
that $K_{G\times G}(X)$ is a module over $K_{T}(\bar{T}^{+})\otimes
R(T)^{W}$ and each filtered piece $F_{\tau}$ is a
$K_{T}(\bar{T}^{+})\otimes R(T)^{W}$-submodule. Furthermore, the
$K_{T}(\bar{T}^{+})\otimes R(T)^{W}$-module structure on $K_{G\times
G}(X)$ given by the above decomposition is compatible with the
canonical $R(T)\otimes R(T)^{W}$-module structure on
$K_{G\times G}(X)$ coming from the inclusion in Prop. \ref{key}.

\ecor {\it Proof:} The existence of the filtration
$\{F_{\tau}\}_{\tau\in\cf_{+}}$ follows by definition. Further, the
filtered pieces multiply by the multiplication rule defined in Theorem
\ref{kdec} and hence it follows that: $F_{\tau}\cdot
F_{\sigma}\subseteq F_{\gamma}$ whenever
$\gamma=\langle\tau,\sigma\rangle$, and $F_{\tau}\cdot F_{\sigma}=0$
whenever $\tau$ and $\sigma$ do not span a cone in $\cf_{+}$.

Recall by Prop. \ref{key} that we have an inclusion
$K_{T}(\bar{T}^{+})\otimes R(T)^{W}\subseteq K_{G\times G}(X)$ as a
subring, which further gives $K_{G\times G}(X)$ a canonical
$R(T)\otimes R(T)^{W}$-module structure.  Now, under the above
isomorphism, $K_{T}(\bar{T}^{+})\otimes R(T)^{W}$ maps to
$\bigoplus_{\tau\in\cf_{+}} C_{\tau}\otimes R(T)^{W}\subseteq
F_{\{0\}}$. Since $F_{\{0\}}\cdot F_{\tau}\subseteq F_{\tau}$, it
follows that each $F_{\tau}$ for $\tau\in\cf_{+}$ is a module over
$K_{T}(\bar{T}^{+})\otimes R(T)^{W}$. Moreover, the decomposition
$K_{G\times G}(X) \simeq \bigoplus_{\tau\in {\cal F}_{+}}
C_{\tau}\otimes R(T)^{W_{\tau}}$, preserves the multiplicative
structure. Thus the above defined $R(T)\otimes R(T)^{W}\subseteq
K_{T}(\bar{T}^{+})\otimes R(T)^{W}$-module structure on $K_{G\times
G}(X)\simeq F_{\{0\}}$, is compatible with the canonical structure
given in Prop. \ref{key}. \hfill $\Box$

The following corollary can be thought of as a geometric
reinterpretation of Theorem \ref{kdec}.

\bcor\label{norm} Let $N_{{\tau}}\simeq\bigoplus_{\rho_j\in\tau(1)}
L_j$ be the normal bundle of $V_{\tau}=\bar{O_{\tau}}$ in $\bar
T$. Let $N_{\tau}\mid_{O_{\tau}}$ denote the restriction of the normal
bundle to $O_{\tau}$ so that
$$\lambda_{-1}(N_{\tau}\mid_{O_{\tau}}):=\prod_{\rho_j\in\tau(1)}
(1-[L_{j}]\mid_{O_{\tau}})\in K_{T}(O_{\tau}).$$ Then we have the
following decomposition:
$$K_{G\times G}(X) \simeq \bigoplus_{\tau\in {\cal F}_{+}}
\lambda_{-1}(N_{\tau}\mid_{O_{\tau}})\cdot K_{T}(O_{\tau}) \otimes
R(T)^{W_{\tau}}.$$ Let
$P_{\tau}:=\lambda_{-1}(N_{\tau}\mid_{O_{\tau}})\cdot K_{T}(O_{\tau})
\otimes R(T)^{W_{\tau}}$ for each $\tau\in {\cal F}_{+}$. Then the
above decomposition is a ring isomorphism where the {\em
multiplication} on the right hand side is given as follows:
\[P_{\tau}\cdot P_{\sigma}\subseteq \left \{\begin{array}{ll}
P_{\gamma}&\mbox{if $\tau$ and $\sigma$ span the cone $\gamma$ in
$\cf_{+}$}\\ 0 &\mbox{if $\tau$ and $\sigma$ do not span a cone in
$\cf_{+}$}\end{array}\right.\] \ecor

{\it Proof:} Observe that $\bz[X_j^{\pm 1}:\rho_j \in \tau(1)]\simeq
K_{T}(O_{\tau})$ since
$K_{T}(O_{\tau})=K_{T}(T/T_{\tau})=R({T}_{\tau})$ where, ${T}_{\tau}$
is the stabilizer of the orbit $O_{\tau}$.  Indeed this isomorphism is
induced from the map (6.2) pp. 27 of \cite{VV} composed with the
restriction to $R({T}_{\tau})$, and is hence compatible with the
isomorphism (2.2) above. Note that under the above isomorphism $X_j$
maps to $[L_j\mid_{O_{\tau}}]$.

Thus $C_{\tau}:=X_{\tau}\cdot \bz[X_j^{\pm 1}:\rho_j \in
\tau(1)]\simeq \lambda_{-1}(N_{\tau}\mid_{O_{\tau}}) \cdot
K_{T}(O_{\tau})=P_{\tau}$, where
$N_{\tau}\mid_{O_{\tau}}\simeq\bigoplus_{\rho_j\in\tau(1)}
L_j\mid_{O_{\tau}}$ denotes the restriction of the normal bundle of
$V_{\tau}$ to $O_{\tau}$.  Therefore, substituting the above
isomorphisms in $(2.3)$ we have the following additive decomposition:

$$ (K_{T}({\bar T}) \otimes R(T))^W \simeq \bigoplus_{\tau\in {\cal
F}_{+}} \lambda_{-1}(N_{\tau}\mid_{O_{\tau}})\cdot K_{T}(O_{\tau})
\otimes R(T)^{W_{\tau}}.$$ Since by Cor.\ref{co2} we have $K_{G\times
G}(X)\simeq (K_{T}({\bar T}) \otimes R(T))^W$, we get the required
decomposition of $K_{G\times G}(X)$.

Further, since the torus embedding $\bar{T}$ is regular, the isotropy
group $T_{\tau}$ has a dense orbit in the normal space $N(x)$ to
$O_{\tau}$ at $x\in O_{\tau}$. In particular, this implies that the
eigenspace of $N(x)$ corresponding to the trivial character of
$T_{\tau}$ is zero. Thus by Lemma 4.2 of \cite{VV} it follows that
$\lambda_{-1}(N_{\tau}\mid_{O_{\tau}})$ is not a zero divisor in
$K_{T}(O_{\tau})$. Thus we see that each piece
$\lambda_{-1}(N_{\tau}\mid_{O_{\tau}})\cdot K_{T}(O_{\tau}) \otimes
R(T)^{W_{\tau}}$ is isomorphic to $K_{T}(O_{\tau}) \otimes
R(T)^{W_{\tau}}$ for every $\tau\in\cf_{+}$.

Furthermore, since $\gamma=\langle \tau,\sigma\rangle$, we have the
restriction maps $R(T_{\gamma})\ra R(T_{\tau})$ $R(T_{\gamma})\ra
R(T_{\sigma})$ induced by the canonical inclusions $T_{\tau}\subseteq
T_{\gamma}$ and $T_{\sigma}\subseteq T_{\gamma}$. These restriction
maps further admit splittings which are respectively given by
$[L_{j}]\mid_{O_{\tau}}\mapsto [L_{j}]\mid_{O_{\gamma}}$ for
$\rho_j\in \tau(1)$ and $[L_{j}]\mid_{O_{\sigma}}\mapsto
[L_{j}]\mid_{O_{\gamma}}$ for $\rho_j\in \sigma(1)$.


In particular,
$\lambda_{-1}(N_{\tau}\mid_{O_{\tau}})=\prod_{\rho_j\in\tau(1)}
(1-[L_{j}]\mid_{O_{\tau}})\in K_{T}(O_{\tau})$ and
$\lambda_{-1}(N_{\sigma}\mid_{O_{\sigma}})=\prod_{\rho_j\in\sigma(1)}
(1-[L_{j}]\mid_{O_{\sigma}})\in K_{T}(O_{\sigma})$ multiply as
elements in $K_{T}(O_{\gamma})$ to give {\footnotesize
$$\prod_{\rho_j\in\tau(1)} (1-[L_{j}]\mid_{O_{\gamma}})\cdot
\prod_{\rho_j\in\sigma(1)} (1-[L_{j}]\mid_{O_{\gamma}})=
\prod_{\rho_j\in\gamma(1)} (1-[L_{j}]\mid_{O_{\gamma}}) \cdot
\prod_{\rho_j\in(\tau\cap \sigma)(1)} (1-[L_{j}]\mid_{O_{\gamma}}).$$}
Thus the right hand side is divisible by
$\lambda_{-1}(N_{\gamma}\mid_{O_{\gamma}})=\prod_{\rho_j\in\gamma(1)}(1-[L_{j}]\mid_{O_{\gamma}})\in
K_{T}(O_{\gamma})$.

Now, by defining multiplication on the right hand side as in Theorem
\ref{kdec} the corollary follows. Thus the above decomposition of
$K_{G\times G}(X)$ is a ring isomorphism. \hfill $\Box$

The structure of rational equivariant cohomology of regular embeddings
has been described in complete detail in \cite{BDP}.  However for
comparison with the setting of $K$-theory, we give below the analogous
statement in the case of cohomology which we obtain by proceeding
along similar steps as in Theorem \ref{kdec}.

We follow the notations in the beginning of this section except for
the following modifications: Let $X_F:=\prod_{\rho_j\in F} X_j$ for
every $F\subseteq \{\rho_j:j=1,\ldots d\}$ in the polynomial algebra
$\bq[X_1,\ldots X_d]$. In particular, let
$X_{\tau}:=X_{\tau(1)}=\prod_{\rho_j\in \tau(1)} X_j$ for every
$\tau\in\cf$.  Let $S:=H^*_T(pt)$ be the symmetric algebra over $\bq$
of $X^*(T)$.  By Theorem 8, pp. 7 of \cite{BDP} we know that:
$$H^*_{T}(\bar T)\simeq \bq[X_1,\ldots X_d]/\langle X_F ~\forall~
F\notin \cf \rangle
$$

Let $e(N_{{\tau}})$ denote the equivariant Euler class of the normal
bundle of $V(\tau)=\bar{O_{\tau}}$ which is equal to the top chern
class of $\bigoplus_{\rho_j\in\tau(1)} L_j$.  We then have the
following (also see Theorem 2.3 of \cite{LP}) for equivariant
cohomology of the wonderful compactification of semisimple adjoint
groups, and the corresponding result for Chow ring of a regular
compactification of a connected reductive group in Cor. 2, p.161 of
\cite{Br}):

\bth{\label{coh}} $$H^*_{G\times G}(X)\simeq (H^*_{T}({\bar T})
\otimes S)^W \simeq \bigoplus_{\tau\in {\cal F}_{+}}
e({N_{\tau}}\mid_{O_{\tau}})\cdot H_{T}^*(O_{\tau}) \otimes
S^{W_{\tau}}
$$ \eeth

\brem As observed in page 3 of \cite{VV} the above results on
algebraic and topological K-theory of $X$ hold with integral
coefficients.  \erem

\subsubsection{Application to Ordinary K-theory}

In this section we shall consider $K_{\tG\times\tG}(X)$, in view of
applying Theorem 4.2 of \cite{Mer} to obtain the results for the
ordinary $K$-ring of $X$. Moreover, by Remark \ref{ext1} we can apply
the contents of \S2 to $K_{\tG\times\tG}(X)$ and
$K_{\tT\times\tT}(X)$.


\bpropo{\label{okth}} Consider the principal ($B^{-}\times B$)-
bundle $G\times G\ra G/B^{-}\times G/B$. Further, we have
a canonical action of $B^{-}\times B$ on $\bar{T}$ through the
surjection $B^{-}\times B\ra T\times T$. We consider the
associated bundle $(G\times G)\times_{B^{-}\times B}{\bar{T}}$
which is a toric bundle with fibre $\bar{T}$ over $G/B^{-}\times
G/B$.  We then have the following description of $K(X)$:

$$K(X)\simeq K((G\times G)\times_{B^{-}\times B}{\bar{T}})^{diag(W)}.$$

\epropo

{\it Proof:} By Cor.\ref{co2} we have: $K_{\tG\times \tG}(X)\simeq
K_{\tT\times \tT}(\bar T)^{diag(W)}.$ Now observe that $K_{\tB^{-}\times
\tB}(\bar T)=K_{\tG\times\tG}((\tG\times\tG)\times_{\tB^{-}\times
\tB}{\bar T})$ (see 5.2.17 of \cite{CG}).  Further, the restriction
homomorphism $K_{\tB^-\times \tB}(\bar{T})\ra K_{\tT\times
\tT}(\bar{T})$ is an isomorphism (see 5.2.18 of \cite{CG}).  Thus we
get the following isomorphism:
$$K_{\tG\times \tG}(X)\simeq
K_{\tG\times\tG}((\tG\times\tG)\times_{\tB^{-}\times \tB}{\bar
T})^{diag(W)}.$$ Note that both sides of the above isomorphism are
$R(\tG)\otimes R(\tG)$-algebras and further, the above isomorphism is
an isomorphism as $R(\tG)\otimes R(\tG)$-algebras (here we use the
fact that $R(\tG)\otimes R(\tG)$ is invariant under the action of
$diag(W)$).  Now, applying the isomorphism $(c)$ of \S1.2 for
$\tG\times \tG$ we get:
$$K(X)\simeq K((\tG\times \tG)\times_{\tB^{-}\times
\tB}{\bar{T}})^{diag(W)} .$$ Further, since the relative $T$-embedding
$(\tG\times \tG)\times_{\tB^{-}\times \tB}{\bar{T}}$, where
$\tB^{-}\times \tB$ acts on the fibre $\bar{T}$ via the surjection to
$B^{-}\times B$, is canonically isomorphic to $(G\times
G)\times_{B^{-}\times B}{\bar{T}}$ we have the proposition.  \hfill
$\Box$

\brem Recall that we have the Cartan decomposition $X=G_{comp}\bar T
G_{comp}$ (see \cite[pp. 585]{DP3}) where $G_{comp}$ is a maximal
compact subgroup of $G$ such that $T_{comp}=T\cap G_{comp}$ is a
maximal compact torus in $G_{comp}$.  Hence for the topological
$K$-theory we have the following isomorphism
$$K^{top}(X)\simeq K^{top}((G_{comp}\times
G_{comp})\times_{T_{comp}\times T_{comp}}{\bar{T}})^{diag(W)},$$ which
is obtained via pullback through the canonical map ${(G_{comp}\times
G_{comp})\times_{T_{comp}\times T_{comp}}\bar{T}}\ra X$ (see \cite
[pp.585-588]{DP3}).  \erem

\brem The above description of $K(X)$ is analogous to the description
of $H^*(X;\bq)$ in Theorem (2.2) of \cite{DP3} in the case when $X$ is
the wonderful compactification. Also see \cite{SU} for the computation
of the Grothendieck ring of a relative torus embedding over an
arbitrary base, analogous to the computation of cohomology in \S3 of
\cite{DP3}. \erem

\brem Here we mention that for the case of a smooth complete toric
variety the structure of the $T$-equivariant and ordinary $K$-theory
is well known (see \S6 of \cite{VV} for the computation of equivariant
and ordinary $K$-theory of any smooth toric variety, and also see
\cite{SU} for the ordinary $K$-theory of a smooth complete toric
variety using different methods).\erem

\section{K-theory of the wonderful compactification}

In this section $X:=\bar{G_{ad}}$ the wonderful compactification of
the semisimple adjoint group $G_{ad}$.

It follows from $(1.1)$ that $G^{ss}$ is the universal cover of
$G_{ad}$, and $T_{ad}:=T^{ss}/C$ is the maximal torus of
$G_{ad}$. Recall that ${rank}~(G_{ad})=rank~(G^{ss})=r$ which is the
semisimple rank of $G$.

The toric variety $\bar{T_{ad}}$ then corresponds to the fan
$\cf_{ad}$ in $X_{*}(T_{ad})\otimes \br$, which is the fan associated
to the Weyl chambers. Moreover, $\bar{T_{ad}}^{+}\simeq\ba^r$, where
$T_{ad}$ acts on $\ba^{r}$ by the embedding $t\mapsto
(t^{\alpha_1},\ldots, t^{\alpha_r})$.  Thus $(\cf_{ad})_{+}$ is the
fan associated to the positive Weyl chamber $\cc^+$ where the edges of
$(\cf_{ad})_{+}$ are generated by the fundamental coweights
$\omega^{\vee}_{1},\ldots, \omega^{\vee}_{r}$ dual to the simple roots
${\alpha_1},\ldots, {\alpha_r}$.

\begin{note}{\label{semisimp}}
Henceforth in this section we let: $G:=G^{ss}$, a semisimple simply
connected algebraic group, $B:=B^{ss}$ a Borel subgroup and
$T:=T^{ss}$ a maximal torus of $G$.
\end{note}

Following Remark \ref{ext1} we shall consider $G\times
G$ and $T\times T$-equivariant $K$-theory of $X$,
where we take the natural actions of $G\times G$ and
$T\times T$ on $X$ through the canonical surjections to
$G_{ad}\times G_{ad}$ and $T_{ad}\times T_{ad}$ respectively. In
particular, we shall apply the contents of \S1 and \S2 for
$K_{G\times G}(X)$ and $K_{T\times T}(X)$.

As in the proof of Prop. \ref{key}, we denote by $u$ and $v$ the first
and second variables of $R(T)\times R(T)$
respectively.

\blem\label{key2} The ring $K_{G\times G}(X)\subseteq
R(T)\otimes R(T)$ consists of elements $f(u,v)\in R(T)\otimes
R(T)$ that satisfy the relations $(1,s_\alpha)f(u,v)\equiv f(u,v)\pmod
{(1-e^{-\alpha(u)})}$ for every $\alpha\in\Delta$.  \elem

{\it Proof:} This follows immediately from the proof of
Prop. \ref{key} since there is only one maximal dimensional cone in
$(\cf_{ad})_+$ which has a facet orthogonal to $\alpha$ for every
$\alpha\in \Delta$. Thus in this case there are no relations of type
$(ii)$ and the relations of type $(i)$ are as given above.\hfill $\Box$

It follows from the above lemma that $R(T)\otimes
R(G)=R(T)\otimes R(T)^{W}\subseteq K_{G\times
G}(X)$ as a subring. In particular, $K_{G\times G}(X)$
is a module over $R(T)\otimes R(G)$, and the following
theorem describes explicitly this module structure.

\bth\label{kwond} The ring $K_{G\times G}(X)$ has the following
direct sum decomposition as $R(T)\otimes R(G)$-module:
$$K_{G\times G}(X)=\bigoplus_{I\subseteq \Delta}\prod_{\alpha\in
I}(1-e^{-{\alpha(u)}}) R(T)\otimes R(T)_{I}.$$ Further, the above
direct sum is a free $R(T)\otimes R(G)$-module of rank $|W|$ with
basis $$\{\prod_{\alpha\in I}(1-e^{-{\alpha(u)}})\otimes f_{v} :~v\in
C^{I}~and~ I\subseteq \Delta\},$$ where $C^{I}$ is as defined in
$(1.3)$ and $\{f_{v}\}$ is as in Notation \ref{aproposds}. \eeth

{\it Proof :} Recall from Lemma \ref{ds} that we have the following
decompositions as $R(T)^W$-modules:
\begin{enumerate}
\item[(i)] $R(T)=\bigoplus_{I}R(T)_{I}$
\item[(ii)] $R(T)^{W_{\Delta\setminus I}}=\bigoplus_{J\subset
I}R(T)_{J}$
\end{enumerate}

Let
$$L:=\bigoplus_{I\subset \Delta}\prod_{\alpha\in I}(1-e^{-{\alpha(u)}})
R(T)\otimes R(T)_{I}$$

For $I\subseteq \Delta$, the piece $\prod_{\alpha\in
I}(1-e^{-{\alpha(u)}}) R(T)\otimes R(T)_{I}$ in the direct sum
decomposition of $L$ is isomorphic to $R(T)\otimes R(T)_I$ as an
$R(T)\otimes R(G)$-module, since $\prod_{\alpha\in
I}(1-e^{-{\alpha(u)}})$ is not a zero divisor in $R(T)$ (see Theorem
\ref{kdec1} for details). Thus it follows from (i) above that $L$ is a
free $R(T)\otimes R(G)$-module of rank $|W|$.

Further, by Lemma \ref{key2} it follows that $\prod_{\alpha\in
I}(1-e^{-{\alpha(u)}}) R(T)\otimes R(T)_{I} \subseteq
K_{G\times G}(X)$ for every $I\subseteq \Delta$. This is
because, $\prod_{\alpha\in I}(1-e^{-{\alpha(u)}})R(T)\otimes
R(T)_{I}\subseteq R(T)\otimes R(T)$ clearly
satisfies the relations which define $K_{G\times G}(X)$
in $R(T)\otimes R(T)$.  In particular, when
$I=\emptyset$, we get $R(T)\otimes R(G)\subseteq
K_{G\times G}(X)$.

Let $K:=K_{G\times G}(X)$. Thus we have the inclusion:
$L\subseteq K$ as modules over $R(T)\otimes R(G)$.

Moreover, by Lemma \ref{fm} we know that $K$ is a free module over
$R(G)\otimes R(G)$ of rank $|W|^2$. Further, note that
$R(T)\otimes R(G)$ is free over $R(G)\otimes R(G)$ of rank
$|W|$. Since $R(G)\otimes R(G)$ and $R(T)\otimes R(G)$ are
regular, it follows that $K$ is a projective module over
$R(T)\otimes R(G)$.  Further, this implies that $K$ is free over
$R(T)\otimes R(G)$ of rank $|W|$, by Theorem 1.1 of \cite{Gu}.

Thus, $L\hookrightarrow K\ra K/L\ra 0$ is a short exact sequence of
$R(T)\otimes R(G)$ modules, and since $K$ and $L$ are
free of rank $|W|$ it follows that $K/L$ is of projective dimension
$1$ as a module over $R(T)\otimes R(G)$.

We require to prove that $L=K$ as $R(T)\otimes
R(G)$-modules. For this we first prove the following lemma.

\blem{\label{loc}} Let $t_{\alpha}:=\prod_{\beta\neq \alpha}
(1-e^{-{\beta(u)}})\in R(T)\otimes R(G)$ for every
$\alpha\in\Delta$. Then, $(K/L)_{t_{\alpha}}=0$ for every
$\alpha\in\Delta$.  \elem

{\it Proof :} Let $M_{\alpha}:=R(T)\otimes R(T)^{s_{\alpha}}
\bigoplus (1-e^{-{\alpha(u)}}) R(T)\otimes
e^{\omega_{\alpha}(v)}.R(T)^{s_{\alpha}}.$

Further note that $R(T)=R(T)^{s_\alpha}\bigoplus
e^{\omega_{\alpha}}R(T)^{s_{\alpha}}$ for every $\alpha\in\Delta$
where $\omega_{\alpha}$ denotes the fundamental weight corresponding
to $\alpha\in\Delta$. Hence, $R(T)\otimes R(T)= R(T)\otimes
R(T)^{s_\alpha}\bigoplus R(T)\otimes
e^{\omega_{\alpha}(v)}R(T)^{s_{\alpha}},$ which is in fact a direct
sum decomposition of $R(T)\otimes R(T)$ as $R(T)\otimes
R(T)^W$-module.

From Lemma \ref{key2} it follows that after localizing at
$t_{\alpha}=\prod_{\beta\neq \alpha} (1-e^{-{\beta(u)}})$, the only
condition defining $K_{G\times G}(X)$ in $R(T)\otimes R(T)$ is
the one corresponding to $\alpha$. Using the above direct sum
decomposition of $R(T)\otimes R(T)$ and the condition
corresponding to $\alpha$, it follows that $K_{t_{\alpha}}\subseteq
(M_{\alpha})_{t_{\alpha}}$.

\noindent
Moreover, from the equalities (i) and (ii) above, we get:\\
$R(T)^{s_{\alpha}}=\bigoplus_{\alpha\notin I} R(T)_{I}$\\
$e^{\omega_{\alpha}}\cdot R(T)^{s_{\alpha}}=\bigoplus_{\alpha\in I}
R(T)_{I}$

\noindent
Hence by the definition of $L$ it further follows that
$L_{t_{\alpha}}= (M_{\alpha})_{t_{\alpha}}.$ Since
$L_{t_{\alpha}}\subseteq K_{t_{\alpha}}\subseteq
(M_{\alpha})_{t_{\alpha}}$, we have: $(K/L)_{t_{\alpha}}=0$ for
every $\alpha\in\Delta$. \hfill $\Box$

Since the projective dimension of $(K/L)=1$, by Auslander-Buchsbaum
formula we know that $Supp~(K/L)$ is of pure codimension $1$ in $Spec~
(R(T)\otimes R(G))$. Hence $Supp~(K/L)$ must contain a prime ideal
${\mathfrak p}$ of height $1$ in $R(T)\otimes R(G)$. Since
$R(T)\otimes R(G)$ is a U.F.D, ${\mathfrak p}=(a)$ for some $a\in
R(T)\otimes R(G)$ and by Lemma \ref{loc} it follows that
${\mathfrak p}$ contains $1-e^{-{\alpha(u)}}$ and $1-e^{-{\beta(u)}}$
for $\alpha\neq\beta\in\Delta$.

This implies that $a$ divides $1-e^{-{\alpha(u)}}$ and
$1-e^{-{\beta(u)}}$ for distinct $\alpha$ and $\beta$, which a
contradiction since $1-e^{-{\alpha(u)}}$ and $1-e^{-{\beta(u)}}$ are
relatively prime in the U.F.D, $R(T)\otimes R(G)$ (see
p.182 of \cite{Bou}).

This contradiction implies that $K/L=0$ and hence $K=L$.

Now, for $I\subseteq \Delta$, the piece $\prod_{\alpha\in
I}(1-e^{-{\alpha(u)}}) R(T)\otimes R(T)_{I}$ in the above direct
sum decomposition is a free $R(T)\otimes R(G)$-module with basis
$\{\prod_{\alpha\in I}(1-e^{-{\alpha(u)}})\otimes f_{v}~:~v\in
C^{I}\}$ where $\{f_{v}\}$ is as in Notation \ref{aproposds}. Thus the
direct sum is a free $R(T)\otimes R(G)$-module of rank $|W|$ with
basis $\{\prod_{\alpha\in I}(1-e^{-{\alpha(u)}})\otimes f_{v} :~v\in
C^{I}~and~ I\subseteq \Delta\}$ \hfill $\Box$

For $I\subseteq \Delta$, let $A_{I}:=\prod_{\alpha\in
I}(1-e^{-{\alpha(u)}}) R(T)\otimes R(T)_{I}\subseteq R(T)\otimes
R(T)$.  The direct sum decomposition in Theorem \ref{kwond} can
therefore be expressed as:
$$K_{G\times G}(X)=\bigoplus_{I\subseteq \Delta}A_{I}.$$

\bcor\label{multwond} The multiplicative structure of $K_{G\times
G}(X)$ is determined by the above decomposition where the pieces
$A_{I}$(resp.  $A_{I^{\prime}}$) corresponding to $I$
(resp. $I^{\prime}$) multiply in $R(T)\otimes R(T)$ as follows:
$$A_{I}\cdot A_{I^{\prime}} \subseteq \prod_{\alpha\in I\cup
I^{\prime}}(1-e^{-{\alpha(u)}}) R(T)\otimes R(T)^{W_{\Delta\setminus
(I\cup I^{\prime})}}\subseteq \bigoplus_{J\subseteq I\cup
I^{\prime}}A_{J}.$$

In particular, any two basis elements $\prod_{\alpha\in
I}(1-e^{-{\alpha(u)}})\otimes f_{v}$ and $\prod_{\alpha\in
I^{\prime}}(1-e^{-{\alpha(u)}})\otimes f_{v^{\prime}}$, where
$v$(resp. $v^{\prime}$) belongs to $C^{I}$ (resp. $C^{I^{\prime}}$)
multiply in $R(T)\otimes R(T)$ to give: {\footnotesize
$$(\prod_{\alpha\in I}(1-e^{-{\alpha(u)}})\otimes f_{v})\cdot
(\prod_{\alpha\in I^{\prime}}(1-e^{-{\alpha(u)}}) \otimes
f_{v^{\prime}})=\prod_{\alpha\in I\cap I^{\prime}}
(1-e^{-{\alpha(u)}})\cdot \prod_{\alpha\in I\cup I^{\prime}}
(1-e^{-{\alpha(u)}})\otimes (f_{v} \cdot f_{v^{\prime}})$$} where the
right hand side belongs to $\prod_{\alpha\in I\cup
I^{\prime}}(1-e^{-{\alpha(u)}}) R(T)\otimes
R(T)^{W_{\Delta\setminus (I\cup I^{\prime})}}$.

\ecor

{\it Proof:} It follows by Lemma \ref{key2} that the direct sum
decomposition $K_{G\times G}(X)=\bigoplus_{I\subseteq
\Delta}A_{I}$ is a ring isomorphism where the multiplication on the
right hand side is given as a subring of $R(T)\otimes R(T)$. We
describe this multiplication below:

Let $B_{I}:=\prod_{\alpha\in I}(1-e^{-{\alpha(u)}}) R(T)\otimes
R(T)^{W_{\Delta\setminus I}}\subseteq R(T)\otimes R(T)$.  By
Lemma \ref{ds} we have $A_{I}\subseteq B_{I}$ and further,
$$B_{I}=\bigoplus_{J\subseteq I} \prod_{\alpha\in I}(1-e^{-{\alpha(u)}})
R(T)\otimes R(T)_{J}\subseteq\bigoplus_{J\subseteq I}
A_{J}.$$ Moreover, we see that $B_{I}$ and $B_{I^{\prime}}$ for $I,
I^{\prime}\subseteq \Delta$ multiply in $R(T)\otimes R(T)$
as follows:
$$B_{I}\cdot B_{I^{\prime}} \subseteq B_{I\cup I^{\prime}}.$$ Thus it
follows that
$$A_{I}\cdot A_{I^{\prime}} \subseteq B_{I\cup I^{\prime}}\subseteq
\bigoplus_{J\subseteq I\cup I^{\prime}}A_{J}.$$
Hence the corollary. \hfill $\Box$

Recall that since $G_{ad}$ is semisimple adjoint,
$\Lambda_{ad}:=X^*(T_{ad})$ has a basis consisting of the simple roots
and further, since $G$ is semisimple and simply connected
$\Lambda:=X^*(T)$ has a basis consisting of the fundamental
dominant weights. Thus $R(T_{ad})=\bz[X^{*}(T_{ad})]$ is generated as
a $\bz$-algebra by $\{e^{\alpha_i}:1\leq i\leq r\}$, and
$R(T)=\bz[X^{*}(T)]$ is generated as a $\bz$-algebra by
$\{e^{\omega_i}:1\leq i\leq r\}$.

Recall that on $X$ the isomorphism classes of line bundles correspond
to $\lambda\in X^*(T)$. Further, the line bundle $\cl_{\lambda}$ on
$X$ associated to $\lambda$, admits a $G\times G$-linearisation so
that $B^{-}\times B$ acts on the fibre $\cl_{\lambda}\mid_{z}$ by
the character $(\lambda,-\lambda)$, where $z$ denotes the base point
of the unique closed orbit $G/B^{-}\times G/B$. Moreover,
$Pic^{G\times G}(X)$ is freely generated by $\cl_{\omega_{i}}$
corresponding to the fundamental dominant weights $\omega_i\in
X^{*}(T)$ for $1\leq i\leq r$ (see \S2.2 of \cite{Br4}).

In particular, $\cl_{\alpha_i}$ are $G\times G$-linearised line
bundles such that $B^{-}\times B$ operates with the character
$(\alpha_i,-\alpha_i)$ on $\cl_{\alpha_i}\mid_{z}$ for $1\leq
i\leq r$. Further, since the centre $\cz$ of $G$ acts trivially on
$X$, and hence acts on the fibre by the character $(\alpha_i,
-\alpha_i)$, $\cl_{\alpha_i}$ is in fact $G_{ad}\times
G_{ad}$-linearised. Moreover, $\cl_{\alpha_i}$ also admits a $G_{ad}\times
G_{ad}$-invariant section $s_i$ whose zero locus is the boundary divisor
$X_i$ for $1\leq i\leq r$.

Moreover, since $\bar{T_{ad}}^{+} \simeq \ba^{r}$ where, each
$I\subseteq \Delta$ corresponds to a $T_{ad}$-orbit
$O_{I}=\{(x_1,\ldots,x_r)\in \ba^r\mid x_i=0$ if and only if
$\alpha_i\in I\}$. The base point $z$ of the closed $G_{ad}\times
G_{ad}$-orbit in $X$ thus corresponds to the $T_{ad}\times
T_{ad}$-fixed point $(0,\ldots, 0)\in \ba^r$. Further, let
$(T_{ad})_I\subseteq T_{ad}$ denote the stabilizer at the orbit
$O_{I}$.

Further, on $\bar{T_{ad}}^{+}$ the line bundle $\cl_{\alpha_i}$ can be
trivialised as a $T_{ad}\times T_{ad}$-equivariant line bundle
$L_{\alpha_i}:=\ba^{r}\times \bc$ where the $T_{ad}\times
T_{ad}$-action is given by $(t_1,t_2)\cdot (u,c)=((t_1,t_2)\cdot u
~,~t_1^{\alpha_i}t_2^{-\alpha_i}\cdot c)$. In particular, we see that
$diag(T_{ad})$ acts trivially on $L_{\alpha_i}$.  The section $s_i$
further becomes the $i$th coordinate function whose zero locus is
$\bar{O_{\{i\}}}$.  Then $N_{I}:=\bigoplus_{\alpha_i\in I}
L_{\alpha_i}$ is the normal bundle of $\bar{O_{I}}$ in
$\bar{T_{ad}}^{+}$, and hence
$$\lambda_{-1}(N_{I}^{\vee}):=\prod_{\alpha_i\in I}
(1-[L_{\alpha_i}^{\vee}]).$$

Note here that $L_{\alpha_i}=L_{i}$ in the notations of \S
2.1. With the above notations we have the following:

\bth\label{norm1} We have the following direct sum decomposition of
 $K_{G\times G}(X)$ as $R(T)\otimes R(G)$-module:
$$K_{G\times G}(X)=\bigoplus_{I\subseteq
\Delta}\lambda_{-1}(N_{I}^{\vee})R(T)\otimes R(T)_{I}. \eqno (3.1)$$
Moreover, $K_{G\times G}(X)$ is free over $R(T)\otimes R(G)$
of rank $|W|$. Further, we can identify the component $R(T)\otimes
1\subseteq R(T)\otimes R(T)^{W}$ in the above direct sum with the
subring of $K_{G\times G}(X)$ generated by
$\{[\cl_{\lambda}]~:~\lambda\in\Lambda\}$, which is also the subring
generated by $Pic^{G\times G}(X)$.  \eeth

{\it Proof:} The isomorphism $(3.1)$ is an immediate consequence of
Theorem \ref{kwond}. We can see this as follows: Since $diag(T)$ acts
trivially on $L_{\alpha_i}$ we see that the isomorphism class of
$L_{\alpha_i}$ in $K_{T\times T}(\bar{T}^{+})$ corresponds to
$e^{\alpha_i}\otimes 1$ in $R(T\times 1)\otimes R(diag(T)$. (Here
we use the canonical identification, $R(T\times \{1\})\otimes
R(diag(T))\simeq R(T)\otimes R(T)$.)  Thus following the
notations in Theorem \ref{kwond}, $[L_{\alpha_i}]$ corresponds to
$e^{\alpha_i(u)}$. Therefore, by the definition of $N_I$ it follows
that the term $\prod_{\alpha\in I}(1-e^{-{\alpha(u)}})$, can be
identified with $\lambda_{-1}(N_{I}^{\vee})$ . Hence the claim.

By Prop. \ref{key} and Remark \ref{ext1} we have the
inclusion $R(T)\times R(G)\subseteq K_{G\times G}(X)$. We
claim that under the above inclusion the image of $R(T)\otimes 1$ in
$K_{G\times G}(X)$ is the subring generated by
$\{[\cl_{\lambda}]~:~\lambda\in\Lambda\}$ which is the subring
generated by the isomorphism classes of $G\times G$-linearised
line bundles on $X$. This can be seen as follows:

By Cor.\ref{co2} and Remark \ref{ext1} we have the canonical inclusion
$K_{G\times G}(X)\hra R(T)\otimes R(T)$ obtained by
restriction to the base point $z$ of the unique closed orbit
$G/B^{-}\times G/B$.  By definition, $\cl_{\lambda}$ maps to
$e^{\lambda}\otimes e^{-\lambda}$ under the above inclusion.  Further,
by the canonical identification $R(T\times \{1\})\otimes
R(diag(T))\simeq R(T)\otimes R(T)$ coming from the exact
sequence $1\ra diag(T)\ra T\times T\ra T\times \{1\}\ra 1$, we
see that the image of $[\cl_{\lambda}]$ under the restriction map is
$e^{\lambda}\otimes 1$ for $\lambda\in \Lambda$.  Since
$e^{\lambda}\otimes 1$ generate $R(T)\otimes 1$, it follows that the
image of $R(T)\otimes 1$ under the above restriction in
$R(T)\otimes R(T)$ is same as the image of the subring generated
by $[\cl_{\lambda}]$ for $\lambda\in\Lambda$.  Further, since
$[\cl_{\lambda}]$ for $\lambda\in\Lambda$ generate
$Pic^{G\times G}(X)$, we have the theorem.  \hfill $\Box$

\brem Note that we can also identify $R(T_{ad})\otimes 1\subseteq
R(T_{ad})\otimes R(T_{ad})^{W}\subseteq K_{G_{ad}\times G_{ad}}(X)$
with the subring of $K_{G_{ad}\times G_{ad}}(X)$ generated by
$\{[\cl_{\alpha_i}]~:~1\leq i\leq r\}$. Thus $K_{G_{ad}\times
G_{ad}}(X)$ is a module over the subring generated by the isomorphism
classes of the $G_{ad}\times G_{ad}$-linearised line bundles on $X$
corresponding to the boundary divisors $\{X_{i}~:~1\leq i\leq
r\}$.\erem

In $R(T)$ let $$f_{v}\cdot f_{v^{\prime}}=\sum_{J\subseteq (I\cup
I^{\prime})} \sum_{w\in C^{J}}a^{w}_{v,v^{\prime}}\cdot f_{w} \eqno
(3.3)$$ for certain elements $a^{w}_{v,v^{\prime}}\in
R(G)=R(T)^{W}$ $~\forall~ v\in C^{I}$, $v^{\prime}\in
C^{I^{\prime}}$ and $w\in C^{J}$, $J\subseteq (I\cup I^{\prime})$ (see
Notation \ref{aproposds}).

\bth\label{kdec1} We have the following isomorphism as $R(T)\otimes
R(T)^{W}$-submodules of $R(T)\otimes R(T)$.
$$ \bigoplus_{I\subseteq \Delta}R(T)\otimes R(T)_{I}\simeq
\bigoplus_{I\subseteq \Delta}\lambda_{-1}(N_{I}^{\vee})R(T)\otimes
R(T)_{I}= K_{G\times G}(X).$$

More explicitly, the above isomorphism maps an arbitrary element
$a\otimes b\in R(T)\otimes R(T)_{I}$ to the element
$(\lambda_{-1}(N_{I}^{\vee})\cdot a)\otimes b$. In particular, the
basis element $1\otimes f_v\in R(T)\otimes
R(T)_{I}$ maps to $\lambda_{-1}(N_{I}^{\vee})\otimes
f_v\in \lambda_{-1}(N_{I}^{\vee})\cdot
R(T)\otimes R(T)_{I}$, for $v\in C^{I}$ for every $I\subseteq
\Delta$.

We now define a multiplication on $\bigoplus_{I\subseteq
\Delta}R(T)\otimes R(T)_{I}$ where any two basis elements
$1\otimes f_v$ and $1\otimes f_{v^{\prime}}$ for $v\in C^{I}$,
$v^{\prime}\in C^{I^{\prime}}$ ($I,I^{\prime}\subseteq \Delta$)
multiply as follows:
$$(1\otimes f_v)\cdot (1\otimes
f_{v^{\prime}}):=\sum_{J\subseteq (I\cup I^{\prime})}\sum_{w\in
C^{J}}(\lambda_{-1}(N_{I\cap I^{\prime}}^{\vee})\cdot
\lambda_{-1}(N_{(I\cup I^{\prime})\setminus J}^{\vee})\otimes
a^{w}_{v,v^{\prime}})\cdot (1\otimes f_{w}).
$$ Then the above isomorphism further preserves the multiplicative
structure where on the right hand side the multiplication is as
defined in Cor.  \ref{multwond}.  \eeth

{\it Proof:} Note that we have a canonical isomorphism of $R(T)$
with $K_{T}(\bar{T}^{+})$ which maps $e^{\alpha_i}$ to
$[L_{\alpha_i}]$ for $1\leq i\leq r$. Thus by the notations in Theorem
\ref{kwond} we have the following identification:
$$\prod_{\alpha\in I} (1-e^{-{\alpha(u)}}) R(T)\otimes
R(T)_{I}=\lambda_{-1}(N_{I}^{\vee})R(T)\otimes R(T)_{I},$$ where
the basis element $\prod_{\alpha\in I}(1-e^{-{\alpha(u)}})\otimes
f_{v}$ corresponds to $\lambda_{-1}(N_{I}^{\vee})\otimes
f_v$ for $v\in C^{I}$ for every $I\subseteq
\Delta$.  Further, in the direct sum decomposition:
$$K_{G\times G}(X)=\bigoplus_{I\subseteq
\Delta}\lambda_{-1}(N_{I}^{\vee})R(T)\otimes R(T)_{I},$$ the
multiplication of two basis elements
$\lambda_{-1}(N_{I}^{\vee})\otimes f_{v}$ and
$\lambda_{-1}(N_{I^{\prime}}^{\vee})\otimes f_{v^{\prime}}$, where
$v$(resp. $v^{\prime}$) belongs to $C^{I}$ (resp. $C^{I^{\prime}}$)
given in Cor.\ref{multwond} can be expressed as follows:
$$(\lambda_{-1}(N_{I}^{\vee})\otimes f_v) \cdot
(\lambda_{-1}(N_{I^{\prime}}^{\vee})\otimes
f_{v^{\prime}})=\lambda_{-1}(N_{I\cap I^{\prime}}^{\vee})\cdot
\lambda_{-1}(N_{I\cup I^{\prime}}^{\vee})\otimes
f_v\cdot f_{v^{\prime}} \eqno(3.4)$$

Note that the right hand side of the above equality $(3.4)$ belongs to
$$\lambda_{-1}(N_{I\cup I^{\prime}}^{\vee}) \cdot R(T)\otimes
R(T)^{W_{\Delta\setminus (I\cup I^{\prime})}}\subseteq
\bigoplus_{J\subseteq (I\cup I^{\prime})}
\lambda_{-1}(N_{J}^{\vee})\cdot R(T)\otimes R(T)_{J},$$ and
further using $(3.3)$ the right hand side of $(3.4)$ can be rewritten
as follows:
$$\sum_{J\subseteq (I\cup I^{\prime})}\sum_{w\in
C^{J}}(\lambda_{-1}(N_{I\cap I^{\prime}}^{\vee})\lambda_{-1}(N_{(I\cup
I^{\prime})\setminus J}^{\vee}) \otimes a^{w}_{v,v^{\prime}})\cdot
(\lambda_{-1}(N_{J}^{\vee})\otimes f_{w}).\eqno(3.5)$$

%

Recall that $\lambda_{-1}(N_{I}^{\vee})$ is not a zero divisor in
$K_{T}(\bar{T_{ad}}^{+})$ (see Lemma 4.2 of \cite{VV} and proof of
Cor. \ref{norm}).  Thus we see that each piece $R(T) \otimes R(T)_{I}$
is isomorphic to $\lambda_{-1}(N_{I}^{\vee})\cdot R(T) \otimes
R(T)_{I}$ for every $I\subseteq \Delta$, as $R(T)\otimes
R(G)$-submodules of $R(T)\otimes R(T)$, where the isomorphism maps an
element $a\otimes b\in R(T) \otimes R(T)_{I}$ to the element
$(\lambda_{-1}(N_{I}^{\vee})\cdot a)\otimes b\in
\lambda_{-1}(N_{I}^{\vee})\cdot R(T) \otimes R(T)_{I}$.

Further, this additively extends to an isomorphism of $R(T)\otimes
R(T)^{W}$-submodules of $R(T)\otimes R(T)$:
$$\bigoplus_{I\subseteq \Delta}R(T)\otimes R(T)_{I}\simeq
\bigoplus_{I\subseteq \Delta}\lambda_{-1}(N_{I}^{\vee})R(T)\otimes
R(T)_{I}.$$
Now by the definition of multiplication in
$\bigoplus_{I\subseteq \Delta}R(T)\otimes R(T)_{I}$ we have:
$$(1\otimes f_v)\cdot (1\otimes
f_{v^{\prime}}):=\sum_{J\subseteq (I\cup I^{\prime})}\sum_{w\in
C^{J}}(\lambda_{-1}(N_{I\cap I^{\prime}}^{\vee})\cdot
\lambda_{-1}(N_{(I\cup I^{\prime})\setminus J}^{\vee})\otimes
a^{w}_{v,v^{\prime}})\cdot (1\otimes f_{w}).
$$
Thus it follows that under the above isomorphism $(1\otimes
f_v) \cdot (1\otimes f_{v^{\prime}})$ maps to
$$\sum_{J\subseteq (I\cup I^{\prime})}\sum_{w\in
C^{J}}(\lambda_{-1}(N_{I\cap I^{\prime}}^{\vee})\lambda_{-1}(N_{(I\cup
I^{\prime})\setminus J}^{\vee}) \otimes a^{w}_{v,v^{\prime}})\cdot
(\lambda_{-1}(N_{J}^{\vee})\otimes f_{w})$$ which by $(3.5)$ is equal
to $(\lambda_{-1}(N_{I}^{\vee})\otimes f_v) \cdot
(\lambda_{-1}(N_{I^{\prime}}^{\vee})\otimes f_{v^{\prime}})$ in
$\bigoplus_{I\subseteq \Delta}\lambda_{-1}(N_{I}^{\vee})R(T)\otimes
R(T)_{I}$. Hence the theorem.\hfill $\Box$

\brem\label{pic} Note that, since
$\lambda_{-1}(N_{\emptyset}^{\vee})=1$, under the isomorphism $$
\bigoplus_{I\subseteq \Delta}R(T)\otimes R(T)_{I}\simeq
\bigoplus_{I\subseteq \Delta}\lambda_{-1}(N_{I}^{\vee})R(T)\otimes
R(T)_{I}$$ defined in Theorem \ref{kdec1} the piece $R(T)\otimes
R(G)=R(T)\otimes R(T)_{\emptyset}$ maps isomorphically onto
itself. In particular, the subring generated by $Pic^{G\times
G}(X)$ in $K_{G\times G}(X)$, which is canonically identified
with $R(T)\otimes \{1\}$ on the right hand side of the above
isomorphism, maps isomorphically onto $R(T)\otimes \{1\}\subseteq
R(T)\otimes R(T)_{\emptyset}$ in the direct sum
$\bigoplus_{I\subseteq \Delta}R(T)\otimes R(T)_{I}$.  \erem

\brem{\label{justkdec1}} Note that Theorem \ref{kdec1} gives an
explicit description of the {\it multiplicative structure constants}
in terms of the basis $\{1\otimes f_{v}:~~v\in C^{I}, I\subseteq
\Delta\}$ for $K_{G\times G}(X)$. This further enables us to {\it
directly apply} this description for the multiplicative structure of
ordinary $K$-ring of $X$ in the following section (see Theorem
\ref{main} and Remark \ref{adveq}).  \erem

\subsection{Ordinary K-ring of the wonderful compactification}

Let $X$ denote the wonderful compactification of $G_{ad}$. We follow
the notations of \S3 (see Notation \ref{semisimp}).

Further, since $K_{G}(G/B)\simeq R(T)$ and $K_{G}(pt)=R(G)$,
the characteristic map $R(T)\ra K(G/B)$ induces an isomorphism
$R(T)/{\cal{J}}\simeq K(G/B)$ where $\cal{J}$ denotes the ideal
generated by $\{f-\epsilon(f)|f\in R(T)^{W}\}$, where
$\epsilon:R(T)\ra \bz$ is the augmentation map given by
$\epsilon(e^{\lambda})=1$ for $\lambda\in \Lambda$.

\blem{\label{pushstein}} Let $\varphi: R(T)\ra K(G/B)$ denote the
characteristic map, and $\bar{f}_v=\varphi(f_v)$ for every $v\in W$.
Then $\bar{f}_v$ for $v\in W$ form a basis of $R(T)/\cj=K(G/B)$ over
$\bz$.  \elem

{\it Proof:} Recall from Notation \ref{aproposds} that $\{f_v\}_{v\in
W}$ form a basis for $R(T)$ as $R(T)^{W}$-module. Since the
characteristic map $\phi:R(T)\ra K(G/B)$ is surjective,
$\{\bar{f}_{v}: v\in W\}$ generate $K(G/B)$ as $\bz$-module. Further,
we claim that $\{\bar{f}_v\}_{v\in W}$ are linearly independent over
$\bz$. This can be seen as follows:

Let $\sum_{v\in W} b_v\cdot \bar{f}_v=0$ for some $b_v\in\bz$. Now,
since $\phi\mid_{R(G)}:R(G)\ra \bz$ is surjective, $b_v=\phi(a_v)$
for some $a_v\in R(G)$ $\forall~v\in W$. Thus we see that
$\sum_{v\in W} a_v\cdot {f}_{v}= c\in \cj$. Further, recall that
$\cj\subseteq R(G)$, and $\{f_{v}:v\in W\}$ are linear independent
over $R(G)$, where $f_1=1\in R(G)$. Thus it follows that $a_1=c$
and $a_v=0$ for all $v\neq 1$. This further implies that $b_v=0$ for
every $v\in W$. Hence the claim. \hfill$\Box$

Further, let
$$K(G/B)_{I}:=\bigoplus_{v\in C^{I}} \bz [{\bar{f}_{v}}],$$
Then we have:
$$K(G/B)=\bigoplus_{I\subseteq \Delta}K(G/B)_{I}.$$

In this section we denote the image in $K(G/B)$ of $e^{\alpha}\in
R(T)$ under the characteristic map by $[L_{\alpha}]$. Further, we
shall denote by the same symbol $\lambda_{-1}(N_{I}^{\vee})\in
K(G/B)$, the image of $\lambda_{-1}(N_{I}^{\vee})=\prod_{\alpha\in
I}(1-e^{-{\alpha}})\in R(T)$ for every $I\subseteq
\Delta$. (However, note that since $\lambda_{-1}(N_{I}^{\vee})\in\cj$,
its image in $K(G/B)$ is always nilpotent.)  Furthermore, we let
$\bar{a}^{w}_{v,v^{\prime}}\in \bz$ denote the image under
$\phi\mid_{R(G)}$ of the element $a^{w}_{v,v^{\prime}}\in
R(G)=R(T)^{W}$ defined in $(3.3)$.

\bth\label{main}We have a canonical $K(G/B)$-module structure on
$K(X)$, induced from the $R(T)\otimes 1$-module structure on
$K_{G\times G}(X)$ given in Theorem \ref{norm1}.  Moreover, $K(X)$
is a free module of rank $|W|$ over $K(G/B)$, $K(G/B)$ being
identified with the subring of $K(X)$ generated by $Pic(X)$.

More explicitly, let
$$\gamma_{v}:=1\otimes \bar{f}_{v}\in K(G/B)\otimes K(G/B)_{I}$$ for
$v\in C^{I}$ for every $I\subseteq \Delta$.  Then we have:

$$K(X)\simeq \bigoplus_{v\in W} K(G/B)\cdot \gamma_v .$$

Further, the above isomorphism is a ring isomorphism, where the
multiplication of any two basis elements $\gamma_v$ and
$\gamma_{v^{\prime}}$ is defined as follows:

$$\gamma_v\cdot \gamma_{v^{\prime}}:=\sum_{J\subseteq (I\cup
I^{\prime})}\sum_{w\in C^{J}}(\lambda_{-1}(N_{I\cap
I^{\prime}}^{\vee})\cdot \lambda_{-1}(N_{(I\cup I^{\prime})\setminus
J}^{\vee})\cdot \bar{a}^{w}_{v,v^{\prime}}) \cdot \gamma_{w}.
$$
\eeth

{\it Proof:} By Theorem \ref{kdec1} we have the following direct sum
decomposition of $K_{G\times G}(X)$ as an $R(T)\otimes
R(G)$-module:

$$K_{G\times G}(X)\simeq\bigoplus_{I\subseteq \Delta}R(T)\otimes
R(T)_{I}.$$

Now, using isomorphism $(c)$ of \S1.2 for $G\times G$ we get:
$$K(X)\simeq \bigoplus_{I\subseteq \Delta} K(G/B)\otimes
K(G/B)_{I}.$$

Further, under the canonical restriction homomorphism $K_{G\times
G}(X)\ra K(X)$, the image of the subring generated by
$Pic^{G\times G}(X)$ in $K_{G\times G}(X)$, maps surjectively
onto the subring generated by $Pic(X)$ in $K(X)$. Hence by Remark
\ref{pic}, it follows that under the above isomorphism, the subring
generated by $Pic(X)$ in $K(X)$ maps isomorphically onto the piece
$K(G/B)\otimes 1\subseteq K(G/B)\otimes K(G/B)_{\emptyset}$.

Let
$$\gamma_{v}:=1\otimes \bar{f}_{v}\in K(G/B)\otimes K(G/B)_{I}$$

for $v\in C^{I}$ for every $I\subseteq \Delta$. Then $\gamma_{v}$ is
the image of the element $1\otimes f_v\in
R(T)\otimes R(T)_{I}$ in $K(G/B)\otimes K(G/B)_I$ under the
characteristic map.

Then identifying $K(G/B)\otimes 1\simeq K(G/B)$, we have:

$$K(X)\simeq \bigoplus_{v\in W} K(G/B)\cdot \gamma_v .\eqno (3.6)$$

Thus we see that $K(X)$ is a free module of rank $|W|$ over $K(G/B)$
with basis $\gamma_v$ for $v\in W$, $K(G/B)$ being identified with the
subring of $K(X)$ generated by $Pic(X)$.

Recall from Theorem \ref{kdec1} that the multiplication of two basis
elements $(1\otimes f_{v})$ and $(1\otimes f_{v^{\prime}})$ of
$K_{G\times G}(X)\simeq\bigoplus_{I\subseteq \Delta}R(T)\otimes
R(T)_{I}$ is defined as:

$$(1\otimes f_v)\cdot (1\otimes f_{v^{\prime}}):=\sum_{J\subseteq
(I\cup I^{\prime})}\sum_{w\in C^{J}}(\lambda_{-1}(N_{I\cap
I^{\prime}}^{\vee})\cdot \lambda_{-1}(N_{(I\cup I^{\prime})\setminus
J}^{\vee})\otimes a^{w}_{v,v^{\prime}})\cdot (1\otimes f_{w}).
$$

Thus their images, $\gamma_v$ and $\gamma_{v^{\prime}}$ in $
\bigoplus_{I\subseteq \Delta} K(G/B)\otimes K(G/B)_{I}$, multiply as
follows:

$$\gamma_v\cdot \gamma_{v^{\prime}}:=\sum_{J\subseteq (I\cup
I^{\prime})}\sum_{w\in C^{J}}(\lambda_{-1}(N_{I\cap
I^{\prime}}^{\vee})\cdot \lambda_{-1}(N_{(I\cup I^{\prime})\setminus
J}^{\vee})\cdot \bar{a}^{w}_{v,v^{\prime}}) \cdot \gamma_{w}.
\eqno(3.7)$$

Thus we conclude that the above isomorphism $(3.6)$ is further a ring
isomorphism, where the multiplication of any two basis elements
$\gamma_v$ and $\gamma_{v^{\prime}}$ for $v\in C^{I}$, $v^{\prime}\in
C^{I^{\prime}}$ and $I,I^{\prime}\subseteq \Delta$ is defined as in
$(3.7)$.

Hence the theorem. \hfill $\Box$

\brem{\label{adveq}} Note that in the direct sum decomposition of
$K_{G\times G}(X)$ given in Theorem \ref{norm}, each piece of the
direct sum is canonically isomorphic to $R(T)\otimes
R(G)$-submodules of $R(T)\otimes R(T)$ (see Theorem
\ref{kdec1}). This enables us in the equivariant setup to describe the
multiplication of the direct sum pieces, and hence the basis elements
inside the subring $R(T)\otimes R(T)$ (see Cor. \ref{multwond} and
Theorem \ref{kdec1}). However, this cannot be done in ordinary
$K$-theory since the image of $\lambda_{-1}(N_{I}^{\vee})$ under the
characteristic homomorphism becomes nilpotent in the ordinary
$K$-ring. Hence the multiplication of the basis elements in ordinary
$K$-ring needs to be defined suitably by pushing down the
multiplicative structure from the equivariant $K$-ring.\erem

\noindent
{\bf Concluding Remarks:}
\begin{enumerate}
\item {\it Extending the results to arbitrary fields and higher
$K$-theory} :

We believe that the results in this paper should hold over any
algebraically closed field of arbitrary characteristic. It is also
likely that many of the results hold in the setting of higher
$K$-theory.

\item {\it Geometric interpretation of the basis $\{f_v\}_{v\in W}$} :

By Prop. \ref{modstein} it follows that when $v\in W^{\dsi}$,
$$\bar{f_{v}}=\sum_{x\in W^{\prime}_{\dsi}(v)}
\varphi(f^{^\emptyset}_{vx})$$ where $\bar{f_{v}}$ and $\varphi$ are
as in Lemma \ref{pushstein}. Further, by (1.5) we see that
$\varphi(f^{^\emptyset}_{vx})$ is the class of the line bundle in
$K(G/B)$ corresponding to the weight
$$v^{-1}\big{(}\sum_{v^{-1}\alpha_i<0}\omega_i\big{)}.$$ We are now
trying to obtain a more comprehensive geometric interpretation of the
basis elements $f_{v}$ and $\bar{f_{v}}$ (and of the Steinberg
basis). Such an interpretation may be well known to experts but we
were unable to find it in the literature.

\item In a private communication we were informed by Prof. De Concini
that E. Strickland has recently determined the structure of the
cohomology ring of the wonderful compactification.

\end{enumerate}

\end{document}